\documentclass[12pt]{article}
\usepackage{amscd}
\usepackage{mathrsfs}
\usepackage{amsfonts}
\usepackage{graphicx}
\usepackage{amsmath,amscd,stmaryrd,amssymb,array, amsfonts,amsmath,amssymb}
\usepackage{enumitem}
\setenumerate[1]{itemsep=0pt,partopsep=0pt,parsep=\parskip,topsep=5pt}
\setitemize[1]{itemsep=0pt,partopsep=0pt,parsep=\parskip,topsep=0pt}
\setdescription{itemsep=0pt,partopsep=0pt,parsep=\parskip,topsep=5pt}

\numberwithin{figure}{section}
 \numberwithin{equation}{section}

\newtheorem{theorem}{Theorem}[section]
\newtheorem{proposition}[theorem]{Proposition}
\newtheorem{definition}[theorem]{Definition}
\newtheorem{corollary}[theorem]{Corollary}
\newtheorem{lemma}[theorem]{Lemma}
\newtheorem{remark}[theorem]{Remark}

\newcommand{\cA}{{\mathcal A}}

\newcommand{\cE}{{\mathcal E}}

\newcommand{\cD}{{\mathcal D}}
\newcommand{\cF}{{\mathcal F}}
\newcommand{\cH}{{\mathcal H}}

\newcommand{\cO}{{\mathcal O}}

\newcommand{\cZ}{{\mathcal Z}}
\newcommand{\cM}{{\mathcal M}}

\newcommand{\cU}{{\mathcal U}}
\newcommand{\cV}{{\mathcal V}}

\newcommand{\cS}{{\mathcal S}}
\newcommand{\cW}{{\mathcal W}}

\newcommand{\cK}{{\mathcal K}}

\newcommand{\sC}{{\mathscr C}}
\newcommand{\sE}{{\mathscr E}}
\newcommand{\sU}{{\mathscr U}}

\newcommand{\sK}{{\mathscr K}}
\newcommand{\sF}{{\mathscr F}}

\newcommand{\mB}{\mb{B}}

\def\be{\begin{equation}}
\def\ee{\end{equation}}
\def\ba{\begin{array}}
\def\ea{\end{array}}
\def\benu{\begin{enumerate}}
\def\eenu{\end{enumerate}}
\def\bt{\begin{theorem}}
\def\et{\end{theorem}}
\def\bp{\begin{proposition}}
\def\ep{\end{proposition}}
\def\bl{\begin{lemma}}
\def\el{\end{lemma}}
\def\br{\begin{remark}}
\def\er{\end{remark}}

\def\b{\beta}
\def\De{\Delta}
\def\de{\delta}
\def\pa{\partial}
\def\nab{\nabla}
\def\lam{\lambda}
\def\Lam{\Lambda}

\def\ve{\varepsilon}
\def\sig{\sigma}

\def\gam{\gamma}

\def\a{\alpha}

\def\.{\cdot}

\def\R{\mathbb{R}}

\def\A{\forall}
\def\ol{\overline}

\def\Cap{\bigcap}\def\Cup{\bigcup}

\def\ra{\rightarrow}

\def\~{\tilde}
\def\8{\infty}
\def\X{\times}
\def\({\left(}
\def\){\right)}

\def\E{\exists}
\def\mb{\mbox}
\def\emp{\emptyset}
\def\sm{\setminus}

\def\Hs{\hspace{1cm}}\def\hs{\hspace{0.5cm}}
\def\Vs{\vskip8pt}\def\vs{\vskip4pt}

\def\({\left(}\def\){\right)}

\begin{document}

\begin{center}
{\bf\Large Local and  Global Dynamic Bifurcations\\[1ex] of Nonlinear Evolution Equations}
\end{center}
\vs\centerline{Desheng  Li
}  
\begin{center}
{\footnotesize
{Department of Mathematics,  Tianjin University\\
          Tianjin 300072,  China\\
{\em E-mail}:  lidsmath@tju.edu.cn}}
\end{center}

\centerline{Zhi-Qiang Wang
\footnote{Corresponding author}
}

\begin{center}
{\footnotesize
{
Center for Applied Mathematics,  Tianjin University\\
     Tianjin 300072,  China, and\\Department of Mathematics and Statistics\\Utah State University, Logan, UT 84322, USA\\
 {\em E-mail}: zhi-qiang.wang@usu.edu}}
\end{center}
\Vs

{\footnotesize
{\bf Abstract.}  We present new local and global dynamic bifurcation results for nonlinear evolution equations of the form $u_t+A u=f_\lam(u)$ on a Banach space $X$, where $A$ is a sectorial operator, and $\lam\in\R$ is the bifurcation parameter.
 Suppose the equation has a trivial solution branch $\{(0,\lam):\,\,\lam\in\R\}$.  Denote  $\Phi_\lam$  the local semiflow generated by the initial value problem of the equation.  It is shown that if the crossing number $n$ at a bifurcation value $\lam=\lam_0$ is nonzero and moreover, $S_0=\{0\}$ is an isolated invariant set of $\Phi_{\lam_0}$, then either there is  a one-sided neighborhood $I_1$ of $\lam_0$  such that $\Phi_\lam$ bifurcates a topological sphere $\mathbb{S}^{n-1}$ for each $\lam\in I_1\setminus\{\lam_0\}$, or there is a two-sided neighborhood $I_2$ of $\lam_0$  such that the system $\Phi_\lam$ bifurcates from the trivial solution  an isolated  nonempty compact invariant set $K_\lam$ with  $0\not\in K_\lam$ for each $\lam\in I_2\setminus\{\lam_0\}$.
We also prove that the bifurcating invariant set has nontrivial Conley index. Building upon this fact we establish a global dynamical bifurcation theorem. Roughly speaking,  we prove that for any given neighborhood $\Omega$ of the bifurcation point $(0,\lambda_0)$, the connected bifurcation branch $\Gamma$ from  $(0,\lambda_0)$ either meets the boundary $\pa\Omega$ of  $\Omega$, or meets another bifurcation point $(0,\lambda_1)$. This result extends the well-known Rabinowitz's Global Bifurcation Theorem to the setting of dynamic bifurcations of evolution equations without requiring the crossing number to be odd.

As an illustration example,  we consider the well-known   Cahn-Hilliard equation.
Some  global  features on dynamical bifurcations of the equation are discussed.

 \Vs

\Vs {\bf 2010 MSC:} 34C23, 34K18, 35B32, 37G99.



}

\newpage

\section{Introduction}

Dynamic bifurcation concerns the changes in the qualitative or topological structures of limiting motions such as equilibria, periodic solutions, homoclinic orbits, heteroclinic orbits and invariant tori etc. for nonlinear evolution equations as some relevant   parameters in the equations vary.
Historically, the subject can be traced back  in the very earlier work of Poincar$\acute{\mb e}$ \cite{Poin}
around 1892. It  is now a fundamental tool to study nonlinear problems in mathematical physics and mechanics \cite{Chow,Kie,Mars}, and enables us to understand how and when a system organizes new states and patterns near the original ``\,trivial\,'' one when the control parameters cross some critical values.

A relatively simpler case for dynamic bifurcation is that of the bifurcations from equilibria.
Generally speaking,   there are two  typical such bifurcations  in the classical bifurcation  theory. One is the bifurcation from equilibria to equilibria (static bifurcation), and the other is from equilibria to periodic solutions (Hopf bifurcation). The former usually requires a ``crossing odd-multiplicity'' condition, namely, the linearized equation of a system has  an odd number of eigenvalues  (counting with multiplicity) crossing the imaginary axis when the control  parameter crosses a critical value (the Krasnosel'skii\,'s Bifurcation Theorem). We also know that in such a case  the bifurcation has some global features, which fact is  addressed by the well-known Rabinowitz's Global Bifurcation Theorem. Situations become very complicated  if one drops  the ``crossing odd-multiplicity'' condition mentioned above. If the system under consideration is a gradient one, then by a classical bifurcation theorem on potential operator equations due to  Krasnosel'skii (see \cite[Chap.\,II, Sect.\,7]{Kie} or \cite{Kras}), one can still have local bifurcation results. 
  Whereas  the global bifurcation remains an open problem. To deal with general systems without the ``crossing odd-multiplicity'' condition,
 Ma and Wang \cite{MW4} proved some new local and global static  bifurcation theorems by using higher-order nondegeneracy conditions on singularities of the nonlinearities. The Hopf bifurcation theory has a long history and,  to some extent, forms the central part of the classical  dynamic bifurcation theory. It focuses on the case when a pair of conjugate eigenvalues  of the linearized equation cross the imaginary axis, and was fully developed in the 20-th century.  There has been  a vast body of literature on how to determine Hopf bifurcation for nonlinear systems arising from applications. One can also find some nice results concerning global results in \cite{AY,Wu}, etc.

This present work is mainly concerned with the general case of the bifurcations from equilibria in terms of invariant-set bifurcation, where the number of  eigenvalues  of the linearized equation crossing the imaginary axis might be even and  greater than   two.
A particular but very important case in this line is the theory of attractor bifurcation, which was first introduced by Ma and Wang in 2003 \cite{MW0}
and was  further developed by the authors into a dynamic transition theory \cite{MW5}.
 Roughly speaking, it states that if the trivial equilibrium solution $\theta$ of a system changes from an attractor to a repeller on the local center manifold  when the bifurcation  parameter $\lam$ crosses a critical value $\lam_0$, then  the system bifurcates a compact invariant set $K$ which is an attractor of the system on the center manifold.
 It is also known that $K$  has the shape of an $n$-dimensional sphere, where $n$ denotes the {\em crossing number} at $\lam=\lam_0$  (the number of  eigenvalues  of the linearized equation crossing the imaginary axis); see \cite[Theorem 1]{san3} or \cite[Theorem 6.1]{MW1}.
  Note  that a fundamental  assumption of this theory is that the trivial equilibrium  $\theta$ is an attractor (repeller) of the system on the center manifold at $\lam=\lam_0$. Hence it is no longer applicable when  $S_0=\{\theta\}$ is only an isolated invariant set when $\lam=\lam_0$.  Fortunately  in such a case, we know that  dynamic  bifurcation still occurs as long as there are eigenvalues crossing the imaginary axis. This  has already been addressed  in the literature; see e.g. Rybakowski \cite{Ryba} (pp.\,101-102) and  Ward \cite{Ward1}.

 An abstract global dynamic bifurcation theorem was also proved in Ward \cite{Ward1} in terms of semiflows on complete metric spaces. Let $\Phi_\lam$ be a family of dynamical systems on a complete metric space $X$, where $\lam\in\R$. Suppose that $\theta$ is an equilibrium solution for each $\Phi_\lam$. Let $[a,b]$ be a compact interval which contains exactly one bifurcation value $\lam_0\in [a,b]$. The Ward's global bifurcation theorem states that  if $h(\Phi_a,\{\theta\})\ne h(\Phi_b,\{\theta\})$, a continua $\Gamma\subset X\X \R$ of bounded solutions bifurcates from $(\theta,\lam_0)$, where  $h(\Phi_\lam,\{\theta\})$ denotes the Conley index of $\{\theta\}$ with respect to $\Phi_\lam$. Moreover, either $\Gamma$ is unbounded in $X\X [a,b]$, or it intersects $X\X\{a,b\}$. 
 Note that, due to the requirement on the uniqueness of  bifurcation values in $[a,b]$, the theorem mentioned above may fail to work  when  a $\lam$-interval contains multiple bifurcation values.  This is somewhat different from the situation of the Rabinowitz's Global Bifurcation Theorem.


In this paper we consider the abstract evolution  equation
\be\label{e:1.1}
u_t+A u=f_\lam(u)
\ee on a Banach space $X$,
where  $A$ is a sectorial operator on $X$ with {\em compact resolvent},   $f_\lam(u)$ is a locally Lipschitz continuous  mapping from $X^\alpha\X \R$ to $X$ for some $0\leq\alpha<1$, and $\lam\in\R$ is the bifurcation parameter. Our main goal is to establish  new local and global dynamic bifurcation results.

Suppose that
\be\label{e:1.2}
f_\lam(0)\equiv0,\Hs \lam\in\R.
\ee
Thus $u=0$ is always a trivial solution of  (\ref{e:1.1}) for all $\lam$.
It is also assumed that $f_\lam(u)$ is differentiable  in $u$ with $Df_\lam(u)$ being continuous in $(u,\lam)$.

First, as one of our main purposes here, we give  some more precise and general results on  local dynamic  bifurcations in terms of invariant sets. In particular, we  show that if the { crossing number} $n$
  at a bifurcation value $\lam=\lam_0$ is nonzero and moreover, $S_0=\{0\}$ is an isolated invariant set of the system, then either  there is  a one-sided neighborhood $I_1$ of $\lam_0$  such that the system bifurcates an $(n-1)$-dimensional  topological sphere $\mathbb{S}^{n-1}$ for each $\lam\in I_1\setminus\{\lam_0\}$, or there is a two-sided neighborhood $I_2$ of $\lam_0$  such that the system bifurcates from the trivial solution  an isolated  nonempty compact invariant set $K_\lam$ with  $0\not\in K_\lam$ for each $\lam\in I_2\setminus\{\lam_0\}$.  

Then  we prove that the invariant set $K_\lam$ from bifurcation has nontrivial Conley index. This result plays a key role in
 establishing our global dynamic  bifurcation theorem. However, it may be of independent interest in its own right.

 Finally, as our main goal in this present work,  we establish a global dynamic bifurcation theorem,   extending the Rabinowitz's Global Bifurcation Theorem on operator equations  to dynamical systems without assuming the ``crossing odd-multiplicity\,'' condition and the uniqueness of bifurcation values in parameter intervals.  Roughly speaking, given a   neighborhood  $\Omega\subset X^\alpha\X\R$ of the bifurcation  point  $(0,\lam_0)$, we  prove that the connected bifurcation branch $\Gamma$ from  $(0,\lam_0)$ either meets the boundary $\pa\Omega$ of $\Omega$, or meets another bifurcation point $(0,\lam_1)$.


As an  example, we consider the homogeneous Neumann boundary value problem of the  Cahn-Hilliard equation
$$u_t+\Delta\(\kappa\De u-f(u)\)=0
$$
 on a bounded domain $\Omega\subset R^d$ ($d\leq 3$) with sufficiently smooth boundary, where
$$f(u)=a_1u+a_2u^2+a_3u^3,\hs \,a_3>0.$$
The local attractor bifurcation and phase transition of the problem have been extensively studied in Ma and Wang  \cite{MW1,MW1b, MW2}. Other results relates to bifurcation of the problem can be found in \cite{BDW,Mis2}, etc.
 Here by applying the  theoretical results obtained above,  we give some more precise  local dynamic bifurcation results and demonstrate  global features of the bifurcations.

 This paper is organized as follows. In Section 2 we make some  preliminaries, and in Section 3 we present results on local invariant manifolds of the equation (\ref{e:1.1}) and give a slightly modified version of a reduction theorem for Conley index in \cite{Ryba}. In Section 4 we prove some  local dynamic bifurcation results.  Section 5 is concerned with the  nontriviality of the Conley indices of bifurcating invariant sets. Section 6 is devoted to the global dynamic bifurcation theorem.  Section 7  consists of  an example mentioned above.


\section{Preliminaries}

This section is concerned with some preliminaries.
\subsection{Basic topological notions and facts}
Let $X$ be  a complete  metric space with metric $d(\cdot,\cdot)$.
For convenience  we will always identify a singleton $\{x\}$ with the
 point $x$ for any $x\in X$.

 Let $A$ and $B$ be nonempty subsets of $X$.
The {\em distance} $d(A,B)$ between $A$ and $B$ is defined as
 $$
 d(A,B)=\inf\{d(x,y):\,\,\,x\in A,\,\,y\in B\},
 $$
 and the {\em Hausdorff semi-distance} and {\em Hausdorff distance} of $A$ and $B$ are defined, respectively,
as
$$d_{\mbox{\tiny H}}(A,B)=\sup_{x \in A}d(x,B),\hs
\delta_{\mbox{\tiny H}}(A,B) = \max\left\{d_{\mbox{\tiny H}}(A,B),
d_{\mbox{\tiny H}}(B,A)\right\}.
$$
 We also assign
$d_{\tiny\mb{H}}(\emp,B)=0$.

 The
closure, interior and boundary of $A$ are denoted, respectively, by $\ol A$, int$\,A$ and $\pa A$. A
subset $U$ of $X$ is called a {\em neighborhood} of $A$, if $\ol A\subset \mbox{int}\,U$. The {\em $\ve$-neighborhood} $\mB(A,\ve)$ of $A$ is defined to be the set $\{y\in
X:\,d(y,A)<\ve\}.$
\vs
Let $A_\lam$ ($\lam\in\Lam$) be a family of nonempty subsets of $X$, where $\Lam$ is a metric space. We say that $A_\lam$ is {\em upper semicontinuous} in $\lam$ at $\lam_0\in\Lam$, this means  $$d_{\tiny\mb{H}}(A_\lam,A_{\lam_0})\ra 0\hs \mb{as }\,\lam\ra\lam_0.$$


\bl\label{l:2.2}\cite{Rab} Let $X$ be a compact metric space, and let $A$ and $B$ be two disjoint closed subsets of $X$. Then either there exists a subcontinuum $C$ of $X$ such that
$$\ba{ll}
A\cap C\ne \emp\ne B\cap C,\ea
$$
or $X=X_A\cup X_B$, where $X_A$ and $X_B$ are disjoint compact subsets of $X$ containing $A$ and $B$, respectively.

\el

\bl\label{l:2.3} (\cite{CV}, pp. 41)\, Let $X$ be a compact metric space. Denote $\sK(X)$ the family of compact subsets of $X$ which is equipped with the Hausdorff metric $\de_{\mbox{\tiny H}}(\.,\.)$. Then $\sK(X)$ is a compact metric space.
\el





\subsection{Criteria  on homotopy equivalence}


We denote ``\,$\simeq$'' and ``\,$\cong$\,''  the {\em homotopy equivalence} and  {\em homeomorphism}, respectively, between topological spaces.


Let $X$ be a topological space, and $A\subset X$ be closed.
The following result  can be found in many text books on general topology.
\vskip-1cm
\bl
If $A$ is a strong deformation retract of $X$, then $X\simeq A.$
\el

Let   $i_A:A\ra X$ be the inclusion. Denote
 $$\ba{ll}
 M_{i_A}=(X\X\{0\})\cup (A\X I),\hs C_{i_A}=M_{i_A}/(A\X\{1\}).\ea
 $$
 $M_{i_A}$ and $C_{i_A}$ are called the {\em mapping cylinder} and {\em mapping cone} of  $i_A$, respectively.

 \vs The pair $(X,A)$ is said to have {\em homotopy extension property} (H.E.P in short), if for any space $Y$, any mapping $f:M_{i_A}\ra Y$ can be extended to a mapping $F:X\X I\ra Y$.
 \bl\label{l2.5} (\cite{Hat}, pp.14)
  $(X,A)$  has the  H.E.P. iff $M_{i_A}$ is a  retract of  $X\X I$.
\el

\bl\label{l2.4}(\cite{Hat}, Theorem 0.17)\, Suppose  $(X,A)$  has the  H.E.P. If $A$ is contractible, then
$
X/A\simeq X.
$
 \el

 As a consequence of Lemma \ref{l2.4}, we have
 \begin{corollary}\label{c2.5} Suppose $(X,A)$  has the  H.E.P. Let $B$ be a closed subset of $A$. If $B$ is a strong deformation retract of $A$,  then
$
X/A\simeq X/B.
$
\end{corollary}
{\bf Proof.} We observe that
$
X/A\cong(X/B)/\widetilde{A},
$ where $\widetilde{A}=\pi_B(A)$, and  $\pi_B:X\ra X/B$ is the projection.
In the following we verify that $$(X/B)/\widetilde{A}\simeq X/B,$$ thus completing the proof of what we desired.

 Since $(X,A)$  has the  H.E.P., $M_{i_A}$ is a  retract of  $X\X I$. Let $f:X\X I\ra M_{i_A}$ be a retraction,
$$
f(x,t)=(\phi(x,t),\xi(x,t)),\Hs (x,t)\in X\X I,
$$
where $\phi(x,t)\in X$, and $\xi(x,t)\in I$.
 Define $$h:X\X I\ra M_{i_{\widetilde{A}}}=((X/B)\X\{0\})\cup (\widetilde{A}\X I)$$ as
$
h(x,t)=\(\pi_B\circ\phi(x,t),\,\xi(x,t)\)$ for $ (x,t)\in X\X I$.
Let $$Q(x,t)=(\pi(x), t),\Hs (x,t)\in X\X I.$$  Then $Q:X\X I\ra (X/B)\X I$ is a quotient mapping.
Observing that
$$
h(x,t)=\(\pi_B\circ\phi(x,t),\xi(x,t)\)= \(\pi_B(x),t\),\Hs (x,t)\in M_{i_A},
$$
one finds that  $h$ remains constant on $B\X \{t\}$ for each $t\in I$. Consequently  $h\equiv const.$  on $Q^{-1}(y,t)$ for each $(y,t)\in (X/B)\X I$. Thus by the basic knowledge in the theory of general topology (see e.g. \cite{Munk}, Chap. 2, Theorem 11.1),
there is a mapping $g:(X/B)\X I\ra M_{i_{\widetilde{A}}}$ such that $h=g\circ Q$. It is trivial to verify that $g$ is  a retraction from $(X/B)\X I$ to $M_{i_{\widetilde{A}}}$.
Thus the pair $(X/B,\widetilde{A})$ has the H.E.P.

 Since $B$  is a strong deformation retract of $A$,  the singleton $\{[B]\}$ is a strong deformation retract of $\widetilde{A}$, that is, $\widetilde{A}$ is contractible.  Lemma \ref{l2.4} then asserts  that $(X/B)/\widetilde{A}\simeq X/B.$ \, $\Box$

\subsection{Wedge/smash product of pointed spaces}
Let $(X,x_0)$ and $(Y,y_0)$ be two pointed  spaces. The {\em wedge product} $(X,x_0)\vee (Y,y_0)$ and {\em smash product} $(X,x_0)\wedge (Y,y_0)$ are defined, respectively, as follows:
$$(X,x_0)\vee (Y,y_0)=\(\cW,\,(x_0,y_0)\),\hs
(X,x_0)\wedge (Y,y_0)=\((X\X Y)/\cW,\,[\cW]\),
$$
where $\cW=X\X\{y_0\}\cup \{x_0\}\X Y$.

We denote  $[(X,x_0)]$ the {\em homotopy type} of a  pointed space $(X,x_0)$.
Since the operations ``$\vee$\,'' and   ``$\wedge$\,'' preserve homotopy equivalence relations, they can be naturally extended to the homotopy types of pointed spaces. Specifically,
$$
[(X,x_0)]\vee [(Y,y_0)]=\left[(X,x_0)\vee(Y,y_0)\right],$$
$$[(X,x_0)]\wedge[(Y,y_0)]=\left[(X,x_0)\wedge(Y,y_0)\,\right].
$$

Denote $\ol0$ and $\Sigma^0$ the homotopy types of the pointed spaces $(\{p\},p)$ and $(\{p,q\},q)$, respectively, where $p$ and $q$ are two distinct points. Let $\Sigma^m$ be the homotopy type of pointed $m$-dimensional sphere. One easily verifies that
$$
[(X,x_0)]\vee \ol0=[(X,x_0)],\hs 
$$
 and
$$
\Sigma^m \wedge\Sigma^n=\Sigma^{m+n},\Hs \A\,m,n\geq 0.
$$

\subsection{Local semiflows and basic dynamical concepts}
 In this subsection we briefly recall some dynamical concepts and facts that will be used throughout the paper.

Let $X$ be a complete metric space.

 A  {\em local semiflow} $\Phi$ on $X$ is a
continuous map from an open subset $\cD_\Phi$ of $\R^+\X X$  to $X$  satisfying that (i)   $\A\,x\in X$,  $\E\,T_x\in(0,\8]$  such that $$(t,x)\in \cD_\Phi \Longleftrightarrow t\in[0,T_x)\,;$$
and (ii) \,$\Phi(0,\.)=id_X$, furthermore,
$$
\Phi(s+t,x)=\Phi\(t,\Phi(s,x)\), \Hs \A\,x\in X,\, \,s,t\geq0
$$ as long as $(s+t,x)\in\cD_\Phi$.
The number $T_x$ in the above definition  is called  the {\em escape time} of $\Phi(t,x)$.
\vs
 Let $\Phi$ be a given local semiflow on $X$. For notational simplicity,  we will rewrite $\Phi(t,x)$ as $\Phi(t)x$.

\vs A {\em trajectory} on an interval $J$  is a continuous  mapping $\gamma: J\ra X$ such that
$$
\gamma(t)=\Phi(t-s)\gamma(s),\Hs \A \,t,\, s\in J,\,\,\,t\geq s.
$$
If $J=\R$, then we simply call  $\gamma$  a {\em complete trajectory}.
The {\em $\omega$-limit set} $\omega(\gamma)$ and {\em $\alpha$-limit set} $\alpha(\gamma)$ of    a { complete trajectory} $\gamma$ are defined, respectively,  as
$$
\omega(\gamma)=\{y:\,\,\,\mb{$\E$ $x_n\in A$ and  $t_n\ra\8$ such that $\gamma(t_n)\ra y$}\},
$$
$$
\a(\gamma)=\{y:\,\,\,\mb{$\E$ $x_n\in A$ and  $t_n\ra-\8$ such that $\gamma(t_n)\ra y$}\}.
$$

Let $S\subset X$. $S$ is said to be   {\em positively invariant} (resp. {\em invariant}), if $\Phi(t)S\subset S$ (resp. $\Phi(t)S=S$) for all $t\geq0$.
A compact invariant set $\cA$ is called an {\em attractor}, if it attracts a neighborhood $U$ of itself, namely,
$$
\lim_{t\ra\8}d_H\(\Phi(t)U,\cA\)=0.
$$
The {\em attraction basin} of an attractor $\cA$, denoted by $\sU(\cA)$, is defined
as
$$\sU(\cA)=\{x:\,\,\lim_{t\ra\8}d(\Phi(t)x,\,\cA)=0\}.$$
\br\label{r2.7}By definition one easily verifies that the {attraction basin} $\sU(\cA)$ of an attractor $\cA$ is open. Furthermore, for any trajectory $\gamma:J\ra X$ of $\Phi$ (where $J$ is an interval), it holds that
$$\mb{ either
$\gamma(J)\subset \sU(\cA)$,\hs or\, $\gamma(J)\cap \sU(\cA)=\emp$}.$$
\er

\subsection{Conley index}
In this subsection we recall  briefly some basic notions and  results in the  Conley index theory. The interested reader is referred to \cite{Conley,Mis} and \cite{Ryba} for details.

Let $\Phi$ be a given local semiflow on $X$, and let $M$ be a subset of $X$.
We say that  {\em $\Phi$ does not explode} in $M$, if  $T_x=\8$ whenever $\Phi([0,T_x))x\subset M.$

  $M$  is said to be  {\em admissible} (see \cite{Ryba}, pp. 13), if for any sequences $x_n\in M$ and $t_n\ra \8$ with $\Phi([0,t_n])x_n\subset M$ for all $n$,  the sequence $\Phi(t_n)x_n$ has a convergent subsequence.
$M$ is said to be {\em strongly admissible}, if it is admissible and moreover, $\Phi$ does not explode in $M$.

\begin{definition}\label{d2.10} $\Phi$ is said to be asymptotically compact on $X$, if each bounded subset $B$ of $X$ is { strongly admissible}.
\end{definition}

From now on   we always assume that
\Vs
({\bf AC}) \,{\em $\Phi$ is asymptotically compact on $X$.}
\Vs
\noindent This  requirement is fulfilled   by a large number of examples from applications.

A compact invariant set $S$ of $\Phi$ is said to be {\em isolated}, if there exists a bounded closed neighborhood $N$ of $S$ such that $S$
is the  maximal invariant set in $N$. Consequently $N$ is called
 an {\em isolating neighborhood} of $S$.

Let there be given an isolated compact invariant set $S$.
 A pair of bounded closed subsets $(N,E)$ is called an {\em  index pair} of $S$, if (i) $N\sm E$ is an isolating neighborhood of $S$;
(ii) $E$ is $N$-invariant, namely, for any $x\in E$ and $t\geq 0$, $$\Phi([0,t])x\subset
N\Longrightarrow\Phi([0,t])x\subset E;$$
(iii) $E$ is an exit set of $N$. That is,  for any $x\in N$, if $\Phi(t_1)x\not\in N$ for some
$t_1>0$, then there exists  $0\leq t_0\leq t_1$ such that
    $\Phi(t_0)x\in E.$

\br 
Index pairs in the terminology of \cite{Ryba} need not be bounded. However, the bounded ones are sufficient for our purposes here.
\er

\begin{definition} (homotopy index) Let $(N,E)$ be an index pair of $S$. Then the
{ homotopy Conley index} of $S$ is defined to be the homotopy
type $[(N/E,[E])]$ of  the pointed space $(N/E,[E])$, denoted by $h(\Phi,S)$.
\end{definition}
\br
Denote  $H_*$ and $H^*$ the singular homology and cohomology theories with coefficient group $\mathbb Z$, respectively.
Applying $H_*$ and $H^*$ to $h(\Phi,S)$ one obtains the {\em homology} and {\em cohomology Conley index} $CH_*(\Phi,S)$ and $CH^*(\Phi,S)$, respectively.
\er

An important property of the Conley index is its continuation property. Here we state a result in this line for the reader's convenience, which is actually a particular case of \cite{Ryba}, Chap.\,\,1,  Theorem 12.2.

 Let $\Phi_\lam$ be a family of semiflows with parameter $\lam\in \Lam$, where $\Lam$ is a connected compact metric space. Assume $\Phi_\lam(t)x$ is continuous in $(t,x,\lam)$.
Denote $\~\Phi$ the {\em skew-product  flow} of the family $\Phi_\lam$ on $X\X \Lam$ defined as follows:
\be\label{e2.3}
\~\Phi(t)(x,\lam)=\(\Phi_\lam(t)x,\lam\),\Hs (x,\lam)\in X\X \Lam.
\ee
\bt\label{t:2.14} Suppose  $\~\Phi$ satisfies the assumption (AC) on $X\X\Lam$. Let $\cS$ be  a compact isolated invariant set of $\~\Phi$.
Then $h(\Phi_\lam,S_\lam)$ is constant for $\lam\in\Lam$, where $S_\lam=\{x:\,\,(x,\lam)\in \cS\}$ is the $\lam$-section of $\cS$.
\et

\noindent{\bf Proof.} Take a bounded closed isolating neighborhood $\cU$ of $\cS$ in $X\X\Lam$. Then the $\lam$-section $\cU_\lam$ of $\cU$ is an isolating neighborhood  of $S_\lam$.  Since $\cS$ is compact in $X\X\Lam$, one easily verifies that $S_\lam$ is upper semicontinuous in $\lam$, namely, $d_{\mbox{\tiny H}}(S_{\lam'},S_\lam)\ra0$ as $\lam'\ra\lam$. Consequently for each fixed $\lam\in\Lam$, $\cU_\lam$ is also an isolating neighborhood  of $S_{\lam'}$ for $\lam'$ near $\lam$. Now the conclusion directly follows from \cite{Ryba}, Chap.\,\,1,  Theorem 12.2.\,  $\Box$

\Vs

Finally, let us also recall the concept of an isolating block.

Let $B\subset X$ be a bounded closed set and $x\in \pa B$ be a boundary point. $x$ is called a {\em strict egress} (resp. {\em strict ingress}, {\em bounce-off}) point of $B$, if for every trajectory $\gamma:[-\tau,s]\ra X$ with $\gamma(0)=x$, where $\tau\geq0$, $s>0$, the following two properties hold.
\begin{enumerate}
\item[(1)] There exists  $0<\ve<s$ such that
$$
\gamma(t)\not\in B \,\,\,(\mb{resp. }\, \gamma(t)\in \mb{int}B,\,\,\,\mb{resp. }\,\gamma(t)\not\in B),\Hs \A\,t\in (0,\ve);
$$
\item[(2)] If $\tau>0$, then there exists $0<\de<\tau$ such that
$$
\gamma(t)\in \mb{int}B \,\,\,(\mb{resp. }\, \gamma(t)\not\in B,\,\,\,\mb{resp. }\,\gamma(t)\not\in B),\Hs \A\,t\in (-\de, t).
$$
\end{enumerate}

Denote  $B^e$ (resp. $B^i$, $B^b$) the set of all strict egress (resp. strict ingress, bounce-off) points of the closed set $B$, and set $B^-=B^e\cup B^b$.

A closed set $B\subset X$ is called an {\em isolating block}   if $B^-$ is closed and $
\pa B=B^i\cup B^-.$
It is well known that if $B$ is a bounded isolating block, then $(B,B^-)$ is an index pair of the maximal compact invariant set $S$ (possibly empty) in $B$.

For convenience,  if $B$ is an isolating block,   we call $B^-$ the {\em boundary exit set}.

\section{Local Invariant Manifolds}

In this section we present some fundamental  results on local invariant manifolds of (\ref{e:1.1}). We also state a slightly modified version of a reduction property of the  Conley index given in \cite{Ryba}.


It is well known that under the hypotheses in Section 1, the initial value problem of (\ref{e:1.1}) is well-posed in $X^\a$. That is, for each  $u_0\in X^\a$ the problem has a unique solution $u(t)$ in $X^\a$ with $u(0)=u_0$ on some maximal existence interval $[0,T)$; see e.g. \cite{Henry}, Theorem 3.3.3.

Denote $\Phi_\lam$ the local semiflow generated by the problem on $X^\a$.
\vs
For convenience in statement, given $Z\subset \mathbb{C}$ and $\a\in\R$, we will write $$\mb{Re}(Z)<\a\,(\,>\a),$$  which  means that
$
\mb{Re}(\mu)< \a\,(\,>\a) $ for all $\mu\in Z.$
\vs



Let $L_{\lam}=A-Df_\lam(0)$. Suppose there exist a neighborhood $J_0=[\lam_0-\eta,\,\lam_0+\eta]$ of $\lam_0\in\R$ and $\de>0$ such that  the following hypotheses  are fulfilled.
\begin{enumerate}
\item[{\bf (H1)}] The spectral $\sig(L_\lam)$  has a decomposition $\ba{ll}\sig(L_\lam)=\Cup_{1\leq i\leq 3}\,\sig_\lam^i\ea $ with
 \be\label{e:2.4}
\mb{Re}(\sig_\lam^1)< -\a_1<-\a_2\leq\mb{Re}\,(\sig_\lam^2) < \a_3<\a_4< \mb{Re}(\sig_\lam^3)\Hs\ee
 for $\lam\in J_0$, where $\a_i$ ($1\leq i\leq 4$) are positive   constants independent of $\lam$.

\vs
\item[{\bf(H2)}]For each $\lam\in J_0$, $X$ has a  decomposition
$X=X^1_\lam\oplus X^2_\lam\oplus X^3_\lam\,$ corresponding to the spectral decomposition in (H1), where
$X^i_\lam$  $(i=1,2,3)$ are $L_{\lam}$-invariant subspaces of $X$. Moreover, $$\mb{dim}\,\,(X^1_\lam),\, \mb{dim}\,(X^2_\lam)<\8.$$

\item[{\bf(H3)}] There is a  family of invertible bounded linear  operators $T=T_\lam$ on $X$  depending continuously on $\lam$ such that when  $\lam\in J_0$, we have
\be\label{e:2.a4}
    T X^i_\lam=X^i_{\lam_0}:=X^i,\Hs i=1,2,3.
  \ee
\end{enumerate}

\br Instead of (H3), a more natural hypothesis  is to assume that\vs
 \begin{enumerate}
\item[{\bf(H3)\,$'$}] the projection operators $P_\lam^{i}:X\ra X^i_\lam$ ($i=1,2$) are continuous in $\lam$.
\eenu\vs
Indeed,  when (H3)\,$'$ is fulfilled, it can be shown that  there is  a  family of invertible bounded linear  operators $T=T_\lam$ on $X$ such that \eqref{e:2.a4} holds true; see \cite{LW2}, Appendix A for details.
\er

We rewrite $E=X^\a$ and set
$$\ba{ll}
E^i=E\cap X^i,\hs E^{ij}=E\cap \(X^i\oplus X^j\),\ea
$$
where  $i,j=1,2,3$ ($i\ne j$).
Then
$$
E=E^2\oplus E^{13}=E^3\oplus E^{12}.
$$
\br
Since $\mb{dim}\,\,(X^1_\lam),\, \mb{dim}\,(X^2_\lam)<\8$, we  have $E^1=X^1,$ $E^2=X^2.$
\er
\bl\label{l:2.1} Assume   (H1)-(H3) are fulfilled.
Then
\begin{enumerate}
\item[(1)]
there exist an open convex neighborhood $W$ of $0$ in $E^2$ and a  mapping $\xi=\xi_{\lam}(w)$ from $W\X J_0$ to $E^{13}$ which is continuous in $(w,\lam)$ and differentiable in $w$, such that for each $\lam\in J_0$,
\be\label{e:2.5}
\cM^2_\lam:=T^{-1}M^2_\lam,\hs \mb{where }\,M^2_\lam:=\{w+\xi_{\lam}(w):\,\,w\in W\},
\ee
is a local invariant manifold of the system (\ref{e:1.1}); and
\item[(2)] there exist an open convex neighborhood $V$ of $0$ in $E^{12}$ and a  mapping $\zeta=\zeta_{\lam}(v)$ from $V\X  J_0$ to $E^3$ which is continuous in $(v,\lam)$ and differentiable in $v$, such that for each $\lam\in J_0$,
\be\label{e:2.5b}
\cM^{12}_\lam:=T^{-1}M^{12}_\lam,\hs \mb{where }\,M^{12}_\lam:=\{v+\zeta_{\lam}(v):\,\,v\in V\},
\ee
is a local invariant manifold of the system (\ref{e:1.1}).
\end{enumerate}
\el

\noindent
{\bf Proof.}  The above results are just slight modifications of the existing ones in the literature; see e.g. \cite{Ryba}, Chap. II, Theorem 2.1.
 Here we give a sketch of the proof for the reader's convenience.

Let $B_\lam=T L_\lam T^{-1}$, and define $$g_\lam(v)=T \(f_\lam(T^{-1} v)-Df_\lam(0)(T^{-1} v)\),\Hs v\in E.$$
Setting $u=T^{-1} v$, the system (\ref{e:1.1}) can be transformed into an equivalent one:
\be\label{e:2.6} v_t+B_\lam v=g_\lam(v).\ee
It is trivial to  check  that $||Dg_\lam(v)||\ra 0$ as $||v||_\a\ra 0$ uniformly with respect to $\lam\in  J_0$. Further by the Mean-value Theorem one easily verifies that for any $\ve>0$, there exists a neighborhood $U$ of $0$ in $E$
such that
\be\label{e:2.7}
||g_\lam(u)-g_\lam(v)||\leq \ve ||u-v||_\alpha,\Hs \A\,u,v\in U,\,\,\lam\in J_0.
\ee

We observe that
$$
B_\lam-\mu I=T L_\lam T^{-1}-\mu I=T (L_\lam-\mu I)T^{-1},
$$
where $I=id_X$ is the identity mapping on $X$, from which it can be easily seen   that $\mu\in \mathbb{C}$ is a regular value of $B_\lam$ if and only if it a regular value of $L_\lam$. Hence  one concludes   that $$\sig(B_\lam)=\sig(L_\lam).$$
Since $X^i_\lam$ ($i=1,2,3$) are $L_{\lam}$-invariant, it follows by (\ref{e:2.a4}) that  $X^i$ are $B_\lam$-invariant for all $\lam\in J_0$.
Now   using some standard argument in the geometric theory of PDEs  (see Henry \cite{Henry}, Sect. 6 and Hale \cite{Hale}, Appendix) and the  uniform contraction principle, it can be shown that
 there exist an open convex neighborhood $W$ of $0$ in $E^2$ and a mapping $\xi=\xi_{\lam}(w)$ from $W\X  J_0$ to $E^{13}$ which is continuous in $(w,\lam)$ and differentiable in $w$, such that for each $\lam\in J_0$,
\be\label{e:3c}
M^2_\lam:=\{w+\xi_{\lam}(w):\,\,w\in W\}
\ee
is a local invariant manifold of the system (\ref{e:2.6}). Consequently $\cM^2_\lam=T^{-1}M^2_\lam$ is a local invariant manifold of  (\ref{e:1.1}).

The proof of the part (2)  follows a fully analogous argument. \, $\Box$
\Vs

  Let $\cM^2_\lam$  and $\cM_\lam^{12}$ be the local invariant  manifolds given in Lemma \ref{l:2.1}, and $\Phi_\lam^2$ and $\Phi_\lam^{12}$ be the restrictions of $\Phi_\lam$ on  $\cM_\lam^2$ and $\cM_\lam^{12}$, respectively, where  $\Phi_\lam$ is the local semiflow  generated by \eqref{e:1.1}.

The following result is  a parameterized  version of  \cite{Ryba}, Chap. II, Theorem 3.1, and can be proved in the same manner as in \cite{Ryba}. We omit the details.

\bl\label{l:3.2}  Assume (H1)-(H3). Then there exist a  neighborhood $U$ of $0$ in $E$ and a number  $\ve>0$ such that for every $\lam\in[\lam_0-\ve,\lam_0+\ve]$,
 \benu
 \item[(1)] $K\subset U$ is a compact  invariant set of $\Phi_\lam$  iff it is a compact invariant set of  $\Phi_\lam^2$ (resp.  $\Phi_\lam^{12}$) on $\cM_\lam^2$ (resp. $\cM_\lam^{12}$); and
 \vs
\item[(2)]  $K\subset U$ is an isolated invariant set of $\Phi_\lam$  iff it is an isolated invariant set of  $\Phi_\lam^2$ (resp.  $\Phi_\lam^{12}$) on $\cM_\lam^2$ (resp. $\cM_\lam^{12}$); furthermore,
  $$
   h\(\Phi_\lam,K \)=  h\(\Phi_\lam^{12},K \)=\Sigma^m\wedge h\(\Phi_\lam^2,K \),$$
where $m=\mb{dim}\,(X^1)$ is the dimension of $X^1$.
\eenu
\el
\vskip-1cm
\section{Local dynamic  bifurcation}
In this section we state and prove some local dynamic bifurcation results concerning   (\ref{e:1.1}) in terms of invariant sets, so we always assume  $$n:=\mb{dim}\,(X^2)\geq 1.$$

In what follows, by a {\em $k$-dimensional topological sphere} we mean the boundary $\pa D$ of any contractible open subset $D$ of a $(k+1)$-dimensional manifold $\cM$ without boundary.
We will  use the notation $\mathbb{S}^k$ to denote any $k$-dimensional topological sphere.

\begin{definition}
 $\mu\in\R$ is called a (dynamic) {bifurcation value} of \eqref{e:1.1}, if  for any neighborhood $U$ of $0$ and $\ve>0$, there exists $\lam\in (\mu-\ve,\mu+\ve)$ such that $\Phi_{\lam}$ has a compact invariant set $K_\lam\subset U$ with $K_\lam\sm\{0\}\ne\emp$.

 If $\mu$ is a bifurcation value, then we call  $(0,\mu)$ a (dynamic) {bifurcation point}.
\end{definition}

We are basically interested in the bifurcation phenomena of the system (\ref{e:1.1}) near a bifurcation value $\lam=\lam_0$. So in addition to  (H1)-(H3), we will also assume
\vs\benu
\item[{\bf(H4)}] there exists $\ve_0>0$ such that for $\lam\in[\lam_0-\ve_0,\,\lam_0+\ve_0]$,
$$
\mb{Re}\,(\sig_\lam^2)<0\,\,(\mb{if }\lam<\lam_0),\hs \mb{Re}\,(\sig_\lam^2)>0\,\, (\mb{if }\lam>\lam_0).
$$
\eenu

Let   $\cM^2_\lam$  and $\cM_\lam^{12}$ be the local invariant  manifolds given in Lemma \ref{l:2.1}, and  $\Phi_\lam^2$ and $\Phi_\lam^{12}$ be the restrictions of $\Phi_\lam$ on  $\cM_\lam^2$ and $\cM_\lam^{12}$, respectively.

\Vs

\noindent{\bf Convention.} \,{\em For simplicity in statement, from now on we  set $\lam_0=0$.}

\subsection{Attractor/repeller bifurcation}


In this subsection we give  an attractor/repeller-bifurcation  theorem, which slightly generalizes some  fundamental  results in Ma and Wang \cite[Theorem 6.1]{MW1} and  \cite[Theorem 4.3]{MW3}.
For the reader's convenience, we also present a self-contained proof for the theorem.
\vs

\bt\label{t:2.2} Assume  (H1)-(H4) are fulfilled (with $\lam_0=0$).

Suppose  $0$ is an attractor (resp. repeller) of $\Phi_0^2$. Then    there exists a closed neighborhood $U$ of $0$ in $E$ and a number $\ve>0$ such that for each
$\lam\in[-\ve,0)$ (resp. $  (0,\ve]$\,), the system $\Phi_\lam$ bifurcates from $0$  a  maximal compact invariant set $K_\lam\ne\emp$ in $U\sm\{0\}$ which contains  an invariant  topological sphere $\mathbb{S}^{n-1}$. Furthermore,
\be\label{e4.5b}\lim_{\lam\ra 0}d_H\(K_\lam,S_0\)=0.\ee

\et


\noindent
{\bf Proof.} {Case 1)} \, $0$ is an attractor of $\Phi_0^2$.
\vs
We first  consider the equivalent system (\ref{e:2.6}) for $\lam\in J_0$.
When  (\ref{e:2.6}) is restricted on  the local center manifold $M^2_\lam$ defined  by (\ref{e:3c}),  it reduces to an ODE system on an open neighborhood $W$ (independent of $\lam$) of $0$ in $E^2$ \,:
\be\label{e:2.9}
w_t=-B_\lam^2 w+P^2\,g_\lam(w+\xi_\lam(w)):=F_\lam(w),\ee
where $B_\lam^2=P^2 B_\lam$, and $P^2$ is the projection from $E=X^\a$ to $E^2$. 
Applying  Lemma \ref{l:3.2} to (\ref{e:2.6}) one deduces that there exist a neighborhood $\cU$ of $0$ in $E$ and $\ve_0>0$ such that for $\lam\in[-\ve_0,\ve_0]$, $S$ is an isolated invariant set of (\ref{e:2.6}) in $\cU$ iff it is an isolated invariant set of the system restricted on the manifold $M^2_\lam$.

Denote $\phi_\lam$ the local semiflow on $W$ generated by (\ref{e:2.9}). Since $0$ is an attractor of $\Phi_0^2$, we find that     $S_0:=\{0\}$ is an attractor of  $\phi_0$. Let $\Omega=\sU(S_0)$  be the attraction basin of $S_0$ in $W$ with respect to $\phi_0$. Then by  converse Lyapunov theorem on attractors (see e.g. \cite[Theorems 3.1 and 3.2]{Li1}),   one can find a function $V\in C^\8(\Omega)$ with $V(0)=0$ and $\lim_{x\ra\pa\Omega}V(x)=+\8$ such that
\be\label{e:2.10}
\nab V(x)\. F_{0}(x)\leq -v(x),\Hs \A x\in\Omega.\ee
where  $v\in C(\Omega)$  and $v(x)>0$ for  $x\ne0$. Let $$N=V_a:=\{x\in\Omega:\,\,V(x)\leq a\}.$$
Then $N$ is a compact  neighborhood of $0$ in $E^2$.
Pick two numbers  $a, \rho>0$ sufficiently small so that \be\label{e4.4}\widetilde{U}:=N\X\mB_{E^{13}}\(\xi_0(N),\rho\)\subset \cU,\ee
where $\xi_0$ is the mapping  determining the local center manifold $M_0^2$ given in Lemma \ref{l:2.1}, and $\mB_{E^{13}}\(\xi_0(N),\rho\)$ denotes the $\rho$-neighborhood of $\xi_0(N)$ in $E^{13}$.

By (\ref{e:2.10}) we have
\be\label{V1}
\nab V(x)\. F_{0}(x)\leq -\mu,\Hs \A x\in\pa N,
\ee
where  $\mu=\min_{x\in \pa N}v(x)>0$, and $\pa N$ is the boundary of $N$ in $E^2$.
Further by  the  continuity of $F_\lam$ in $\lam$,  there exists $0<\ve_1\leq \ve_0$ such that
\be\label{e:3.b}
\nab V(x)\. F_{\lam}(x)\leq -\frac{\mu}{2},\Hs \A x\in\pa N
\ee
for $\lam\in[-\ve_1,\ve_1]$, which implies   that $N$ is a positively invariant set of  $\phi_\lam$.

 It can be  assumed that  $\ve_1$ is  sufficiently small  so that
\be\label{e:3.c}
\xi_\lam(N)\subset\mB_{E^{13}}\(\xi_0(N),\rho\),\Hs \lam\in[-\ve_1,\ve_1].
\ee

Now assume  $\lam\in[-\ve_1,0)$.   Consider   the inverse flow $\phi_\lam^-$ of $\phi_\lam$  generated by the system
\be\label{e:2.11}
w_t=-F_\lam(w):=B_\lam^2 w-P^2\,g_\lam(w+\xi_\lam(w)).\ee
By (H4) we find that $\mb{Re\,}(\sig(B_\lam^2))<0,$  which implies that    $S_0$ is an attractor of  $\phi_\lam^-$.
  Let $G_\lam=\sU(S_0)$ be the attraction basin of $S_0$ in $W$ with respect to  $\phi_\lam^{-}$.
  We infer from (\ref{e:3.b}) that each  $x\in \pa N$ is a strict ingress point of $\phi_\lam$, and hence is a strict egress point of $\phi_\lam^-$. Thus  one necessarily has $G_\lam\subset N$. Therefore the boundary $\pa G_\lam$ of $G_\lam$ in $E^2$ is contained in $N$; see Fig. \ref{fg3-1}.

  \begin{figure}[h]\centering\includegraphics[width=8.5cm]{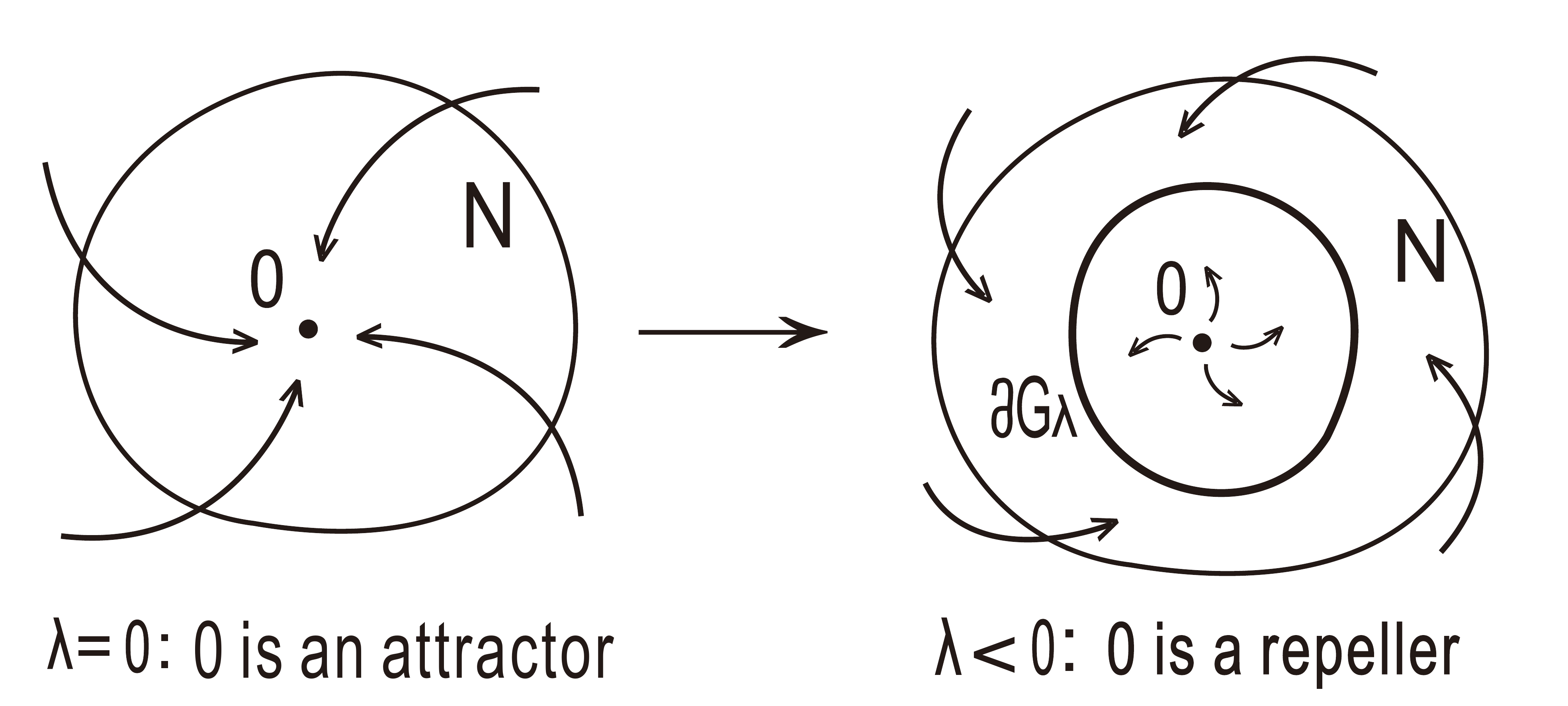}\vskip-0.3cm \caption{Attractor-bifurcation}\label{fg3-1}\end{figure}

  We prove that $\pa G_\lam$   is an invariant set of $\phi_\lam^-$.  For this purpose,  it suffices to show that for each $x_0\in \pa G_\lam$, there is a complete trajectory $w(t)$ of $\phi_\lam^-$ (\,i.e., a solution of \eqref{e:2.11}\,) with $w(0)=x_0$ such that $w(t)\in \pa G_\lam$ for all $t\in\R$.

Note that \eqref{e:2.11} always has a unique   solution $w(t)$   defined on a maximal existence interval $J$ such that $w(0)=x_0$. Since $x_0\not\in G_\lam$, by Remark \ref{r2.7} we deduce that $w(t)\not\in G_\lam$ for all $t\in J$. We claim that $w(t)\in \pa G_\lam$ for $t\in J$, and consequently one also has $J=\R$, thus completing  the proof of the invariance of $\pa G_\lam$.
  We argue by contradiction and suppose the claim was false.  Then there would exist $t_0\in J$ such that $w(t_0)\not\in \ol{G}_\lam$. Hence $d\(w(t_0),\ol{G}_\lam\)>0$. Take a sequence $x_k\in G_\lam$ such that $x_k\ra x_0$.
  Let $w_k(t)$ by the solution of \eqref{e:2.11} with $w_k(0)=x_k$. Then by continuity properties on ODEs, we know that $t_0$ belongs to the maximal existence interval $J_k$ of $w_k(t)$ if $k$ is sufficiently large; furthermore, $w_k(t_0)\not\in \ol{G}_\lam$.
But  by Remark \ref{r2.7},   this leads to a contradiction since $w_k(0)=x_k\in G_\lam$.

  Denote $A_\lam$ the maximal compact invariant set of $\phi_\lam$ in $N\sm G_\lam$. Clearly  $\pa G_\lam\subset A_\lam$. It is trivial to check that $A_\lam$ is  the maximal compact invariant set of $\phi_\lam$ in $N\sm S_0$.
Since $N$ is an isolating neighborhood of $S_0$ with respect to $\phi_0$, by a simple argument via contradiction it can be shown that
\be\label{e4.5a}\lim_{\lam\ra 0}d_H\(A_\lam,S_0\)=0.\ee

 We claim that  $\pa G_\lam$ is an $(n-1)$-dimensional topological sphere. Indeed, define
$$
H(s,x)=\left\{\ba{ll}\phi_\lam^{-}\(\frac{s}{1-s}\)x,\Hs\, &s\in[0,1),\,\,x\in G_\lam;\\[1ex]
0,& s=1,\,\,x\in G_\lam.\ea\right.
$$
Then $H$ is a strong deformation retraction shrinking $G_\lam$ to the point $0$. This shows that $G_\lam$ is contractible and proves our claim.

Now we define $$\widetilde{K}_\lam=\{w+\xi_\lam(w):\,\,w\in A_\lam\},\hs \widetilde{\mathbb{S}}=\{w+\xi_\lam(w):\,\,w\in\pa G_\lam\},$$
where $\xi_\lam$ is the mapping  in \eqref{e:2.9} given by Lemma \ref{l:2.1}.
By (\ref{e:3.c}) and  (\ref{e4.4}) we find that $\widetilde{K}_\lam\subset \widetilde{U}\subset \cU$.  $\widetilde{K}_\lam$ is the maximal compact invariant set of the system (\ref{e:2.6}) in $\widetilde{U}\setminus \{0\}$.
It follows by \eqref{e4.5a}   that $\lim_{\lam\ra 0}d_H(\widetilde{K}_\lam,S_0)=0$.
\vs
Finally, let
  $U_\lam=T^{-1}\widetilde{U},$ where $T=T_\lam$ is the linear operator in (H3). Then one can find a closed neighborhood $U$ of $0$ in $E$ and a number  $0<\ve\leq\ve_1$ such that $U\subset U_\lam$ for all $\lam\in[-\ve,0)$.
    Set $ K_\lam=T^{-1}\widetilde{K}_\lam$. Then $\lim_{\lam\ra 0}d_H({K}_\lam,S_0)=0$. Thus we may assume that $\ve$ is chosen sufficiently small so that $K_\lam\subset \mb{int}\,U$ for all $\lam\in[-\ve,0)$.
 It is easy to see that  $U$ and $K_\lam$    fulfill all the requirements of the theorem.

\vs
{\bf Case 2)} \,The equilibrium $0$ is a repeller of $\Phi_0^2$.

\vs
This case  can be treated by  replacing  (\ref{e:2.9}) and (\ref{e:2.11}) with each other and repeating the  above argument. We omit the details.\, $\Box$

\subsection{Invariant-set bifurcation}

Now we state and prove a  general local invariant-set bifurcation theorem.

\bt\label{t:2.7} Assume that  (H1)-(H4) are fulfilled. Suppose $S_0=\{0\}$ is an  isolated invariant set of $\Phi_0$.
 Then    one of the following assertions holds.
\begin{enumerate}
\item[(1)]  $S_0$ is an attractor (resp. repeller) of $\Phi_0^2$. In such a case, the system undergoes an attractor-bifurcation (resp. repeller-bifurcation)  in Theorem \ref{t:2.2}.
\vs
    \item[(2)] There exist a closed  neighborhood $U$ of $0$ in $E$ and a two-sided neighborhood $I_2$ of $\lam_0$  such that $\Phi_\lam$ has a nonempty maximal compact invariant set $K_\lam$ in $U\sm S_0$ for each $\lam\in I_2\sm\{\lam_0\}$.
\end{enumerate}
 Furthermore,  in both cases the bifurcating invariant set $K_\lam$  is upper semicontinuous  in ${\lam}$ with  \be\label{e3.18a}\lim_{\lam\ra0}d_{\mbox{\tiny H}}\(K_\lam,\,0\)=0.\ee
\et


\noindent
{\bf Proof.} Let us first verify the bifurcation results  in (1) and (2). For this purpose, it suffices to assume  $S_0$ is neither an attractor nor a repeller of $\Phi_0^2$ and prove that the second assertion (2) holds true.
\vs

Let us start with   the local semiflow $\phi_\lam$ generated  by the bifurcation equation  (\ref{e:2.9}) on $W$.
Since $S_0=\{0\}$ is an  isolated invariant set of $\Phi_0$, by Lemma \ref{l:3.2} it is isolated for $\Phi_0^2$.  Because $\Phi_\lam^2$  and   $\phi_\lam$  are   conjugate, one concludes that $S_0$ is an  isolated invariant set of $\phi_0$.

Note that (H4) implies  $$\mb{Re\,}(\sig(B_\lam^2))<0\,\,\,(\lam<0),\hs \mb{Re\,}(\sig(B_\lam^2))>0\,\,\,(\lam>0),$$
where $B_\lam^2$ is the linear operator in \eqref{e:2.9}. Hence  $S_0$ is  a repeller of $\phi_\lam$ when $\lam<0$, and an attractor when $\lam>0$.
By  Lemma \ref{l:3.2} we also have for some $\ve_1>0$ that
\be\label{e:4.5}
\Hs h\(\phi_\lam,S_0\)=\Sigma^{n}\,\,\,(\lam\in[-\ve_1,0)),\hs h\(\phi_\lam,S_0\)=\Sigma^{0}\,\,\,( \lam\in(0,\ve_1]).
\ee

Pick a closed  neighborhood $W_0$ of $S_0$ in $E^2$ such that it is an isolating neighborhood of  $S_0$ with respect to $\phi_0$.  Then by a simple argument via contradiction,  we deduce that  $W_0$ is also an isolating neighborhood of the maximal compact invariant set $S_\lam$ of $\phi_\lam$ in $W_0$ provided $\lam$ is sufficiently small; furthermore, \be\label{e3.18}  \lim_{\lam\ra0}d_{\mbox{\tiny H}}\(S_\lam,\,S_0\)=0.  \ee

Fix a positive number $\ve<\ve_1$ such that $W_0$ is  an isolating neighborhood  of  $S_\lam$ for all $\lam\in [-\ve,\ve]$.
Then Theorem \ref{t:2.14} asserts that
\be\label{e:2.13}
h(\phi_\lam,S_\lam)\equiv const.,\Hs \lam\in[-\ve,\ve].
\ee

In what follows we show that  \be\label{e:4.a1}h\(\phi_0,S_0\)\ne\Sigma^{0}.\ee

Since $S_0$ is an isolated  invariant set of $\phi_0$, by \cite{CE}, Theorem 1.5, one can find  a connected  isolating block $B$ of $S_0$ (with respect to $\phi_0$) with smooth boundary $\pa B$.
We claim that  $B^-\ne \emp$, where $B^-$ is the boundary exit set of $B$ with respect to the flow $\phi_0$. Indeed, if $B^-= \emp$ then $B$ is positively invariant under the system $\phi_0$. Because $S_0$ is the maximal compact invariant set of $\phi_0$ in $B$, one easily deduces that it is an attractor of $\phi_0$, which contradicts the assumption that $S_0$ is not an attractor of $\Phi_0^2$ (recall that $\Phi_\lam^2$ and $\phi_\lam$ are conjugate).

Denote  $H_*$  the singular homology theories with coefficient group $\mathbb Z$. Then  $h\(\phi_0,S_0\)=[(B/B^-,[B^-])]$. Therefore
$$
H_0(h\(\phi_0,S_0\))=H_0((B/B^-,[B^-]))=H_0(B,B^-).
$$
As $B$ is path-connected and $B^-\ne\emp$, by the basic knowledge in the theory of algebraic topology we find that $H_0(B,B^-)=0$. Consequently
$
H_0(h\(\phi_0,S_0\))=0$. On the other hand, recalling that $\Sigma^{0}$ is the homotopy type of any pointed space $(\{p,q\},q)$ consisting of exactly two distinct points  $p$ and $q$, one has
$$
H_0(\Sigma^{0})=H_0((\{p,q\},q))=\mathbb{Z}.
$$
Hence we see that \eqref{e:4.a1} holds true.

\vs Now assume  $\lam\in(0,\ve]$.  Combining  (\ref{e:4.5}) and  (\ref{e:2.13}) it yields
$$
h(\phi_\lam,S_\lam)=h(\phi_0,S_0)\ne h\(\phi_\lam,S_0\),$$
which implies  that
$S_\lam\sm S_0\ne\emp.$ Recall that  $S_0$ is an attractor of $\phi_\lam$. Let $$R_\lam=\{x\in S_\lam:\,\,\omega(x)\cap S_0=\emp\}.$$ Then  $R_\lam$ is a nonempty compact invariant set of  $\phi_\lam$ with $(R_\lam,S_0)$ being a repeller-attractor pair of $S_\lam$; see \cite{Ryba},  pp.141. Because  $S_\lam$ is maximal in $W_0$, it can be easily seen that  $R_\lam$ is precisely the maximal compact invariant set in $W_0\sm S_0$.

Consider the inverse flow $\phi_\lam^-$ of $\phi_\lam$ on $W$. Then we have \be\label{e:4.5b}
\hs h\(\phi_\lam^-,S_0\)=\Sigma^{0}\,\,\,(\lam\in[-\ve,0)),\hs h\(\phi_\lam^-,S_0\)=\Sigma^{n}\,\,\,( \lam\in(0,\ve]).
\ee
Since $S_0$ is a repeller of $\phi_\lam$ for $\lam\in[-\ve,0)$, it is an attractor of $\phi_\lam^-$. Repeating the same argument above with $\phi_\lam$ replaced by $\phi_\lam^-$, one immediately deduces that $\phi_\lam$ has a nonempty maximal compact invariant set $R_\lam$  in $W_0\sm S_0$ for $\lam\in[-\ve,0)$.

\vs

We show that $R_\lam$ is upper semicontinuous  in $\lam$.
We only consider the case where
$\lam\in(0,\ve]$. The argument for the case where $\lam\in[-\ve,0)$ can be performed in the same manner by considering the inverse flow $\phi_\lam^-$, and we omit the details.

 Let $\sU_\lam=\sU(S_0)$ be  the attraction basin of $S_0$ in $W$ with respect to $\phi_\lam$.
 For each fixed $\lam> 0$, pick a number  $r>0$ such that  $\ol\mB_r\subset \sU_\lam$, where (and below)  $\mB_r$  denotes the ball in $E^2$ centered at $0$ with radius $r$. Then by the stability property of  attraction basins (see e.g. Li \cite[Theorem 2.9]{Li0}), there exists $\rho>0$ such that  $\ol\mB_{r/2}\subset \sU_{\lam'}$ provided $|\lam'-\lam|\leq\rho$.
 This implies that  $$\ba{ll}R_{\lam'}\cap \ol\mB_{r/2}=\emp\ea$$ for all $\lam'\in(0,\ve]$ with $|\lam'-\lam|\leq\rho$. We check that $
\lim_{\lam'\ra\lam}d_{\mbox{\tiny H}}\(R_{\lam'},\,R_\lam\)=0,
$
thus proving  what we desired.

Suppose the contrary. There would exist $\lam_k\ra\lam$ and $\de_0>0$ such that $$d_{\mbox{\tiny H}}\(R_{\lam_k},\,R_\lam\)\geq \de_0,\Hs\A\,k\geq 1.$$ We may assume $|\lam_k-\lam|\leq\rho$, hence $R_{\lam_k}\subset W_0\sm\mB_{r/2}$ for all $k$. Thanks to Lemma \ref{l:2.3}, it can be assumed that $R_{\lam_k}$  converges to a nonempty compact subset $R'_\lam$ of $W_0\sm\mB_{r/2}$ in the sense of Hausdorff distance $\de_{\mbox{\tiny H}}(\.,\.)$. Then $d_{\mbox{\tiny H}}\(R'_\lam,\,R_\lam\)\geq \de_0$. On the other hand, one trivially verifies that $R'_\lam$ is an invariant set of $\phi_\lam$. Thus $R_\lam\cup R'_\lam$ is a compact invariant set of $\phi_\lam$ in $W_0\sm S_0$. This contradicts  the maximality of $R_\lam$ in $W_0\sm S_0$.

We are now ready to complete the proof of the theorem. Let  $U$ be the neighborhood of $0$ given in Lemma \ref{l:3.2}. We may restrict $U$ sufficiently small in advance so that $P^2T_\lam U\subset W_0$ for all $\lam\in[-\ve,\ve]$, where $P^2:E\ra E^2$ is the projection, and $T_\lam$ is the operator in (H3). Let $K_\lam=T_\lam^{-1}\widetilde{R}_\lam$, where $$\widetilde{R}_\lam=\{w+\xi_\lam(w):\,\,w\in R_\lam\}.
$$
Then $K_\lam$ is upper semicontinuous in $\lam$ and  is a  compact invariant set of $\Phi_\lam$. By \eqref{e3.18}  we have $\lim_{\lam\ra0}d_{\mbox{\tiny H}}\(R_\lam,\,S_0\)=0$. It follows that $\lim_{\lam\ra0}d_{\mbox{\tiny H}}(K_\lam,\,S_0)=0$. Thus one can assume  $\ve$ is chosen sufficiently small  so that $K_\lam\subset U$ for  $\lam\in[-\ve,\ve]$.

We claim that $K_\lam$ is the maximal compact invariant set of $\Phi_\lam$ in $U\sm S_0$, which completes the proof of the theorem. Indeed, if this was false,
then $\Phi_\lam$ would have another compact invariant set $K'_\lam\subset U\sm S_0$ such that $K_\lam \varsubsetneq K'_\lam$. It follows  that
$$
R_\lam=P^2T_\lam K_\lam \varsubsetneq P^2T_\lam K'_\lam:=R'_\lam.
$$
By the invariance of $K'_\lam$ it is easy to deduce that $R'_\lam$ is a compact invariant set of $\phi_\lam$ in $W_0\sm S_0$. However, this contradicts the maximality of $R_\lam$ in $W_0\sm S_0$.
 \,$\Box$
\subsection{Some remarks on static bifurcation}

It is worth noticing that Theorem \ref{t:2.7} may also give us  information on
the static  bifurcation of the system in some cases.
For example, if the stationary problem
\be\label{3.27}
Au=f_\lam(u), \Hs u\in E:=X^\alpha
\ee
  has a variational structure, then  (\ref{e:1.1}) is a gradient-like system, and each nonempty compact invariant set $K$ of $\Phi_\lam$ contains at least one equilibrium point, which is precisely a solution of (\ref{3.27}). On the other hand, it is also easy to see that  if $K$ consists of at least two distinct points, then it contains at least two distinct equilibrium points of $\Phi_\lam$. Thus under the hypotheses of Theorem \ref{t:2.7}, one immediately concludes  that either   there is a one-sided neighborhood $I_1$ of $\lam_0$ such that (\ref{3.27}) bifurcates two distinct nontrivial solutions for each $\lam\in I_1\sm\{\lam_0\}$, or there is a two-sided neighborhood $I_2$ of $\lam_0$ such that (\ref{3.27}) bifurcates at least one nontrivial solution for each $\lam\in I_2\sm\{\lam_0\}$.

We refer the interested reader to \cite{CW,Rab2,RSW} and \cite{SW}, etc. for more detailed  bifurcation results on such operator equations.

\vs
 As another  example, we consider the  particular but important case  where $$n=\mb{dim}(X^2)=1.$$ We first claim that each compact invariant set $C_\lam$ of $\Phi_\lam$ close to $0$ contains at least one equilibrium point which is a solution of (\ref{3.27}). Indeed, each such invariant set $C_\lam$  is contained in the local invariant  manifold $\cM_\lam^2$. Because  $\cM_\lam^2$ is a $C^1$ curve,  every connected component $\ell$ of  $C_\lam$ is a segment of $\cM_\lam^2$. Since (\ref{e:1.1}) reduces to a one-dimensional ODE on $\cM_\lam^2$ (hence backward uniqueness holds on $\cM_\lam^2$),
 by invariance of $\ell$ it is trivial to deduce that   the end points of $\ell$ are   equilibria of $\Phi_\lam$.

Using the above basic fact, we can also easily verify that  $0$ is an isolated solution of (\ref{3.27}) at $\lam_0$ if and only if  $S_0=\{0\}$ is an isolated invariant set of $\Phi_{\lam_0}$.
Thanks to Theorem \ref{t:2.7}, one immediately obtains  the following bifurcation result, which generalizes Henry \cite{Henry}, Theorem 6.3.2.

\bt\label{t:3.6} Assume  (H1)-(H4) are fulfilled with $\mb{dim}\,(X^2)=1$.  Then one of the  following alternatives occurs.
\begin{enumerate}

\item[(1)] There is a sequence  $u_k$ of nontrivial solutions  of (\ref{3.27}) at $\lam=\lam_0$ such that $u_k\ra 0$ as $k\ra\8$.

\item[(2)] There is a one-sided neighborhood $I_1$ of $\lam_0$ such that  (\ref{3.27}) bifurcates  at least two nontrivial solutions for each  $\lam\in I_1\sm\{\lam_0\}$.

\item[(3)] There is a two-sided neighborhood $I_2$ of $\lam_0$ such that  (\ref{3.27}) bifurcates  at least one nontrivial solution for each  $\lam\in I_2\sm\{\lam_0\}$.
\end{enumerate}

\et




\br
When $\mb{dim}\,(X^2)=1$ we can also use the classical Crandall-Rabinowitz Theorem (see \cite{Kie}, Theorem I.5.1)
to derive more explicit static  bifurcation results  under some additional assumptions such as the {\em transversality condition}.
(Some nice bifurcation results when the transversality condition mentioned above is violated can be found in \cite{Shi} etc.) Other general bifurcation theorems such as the Krasnosel'skii Bifurcation Theorem (see \cite{Kie}, Theorem II.3.2) also apply to deal with this special case.

\er

\br
Whether the bifurcating  invariant set $K_\lam$   contains equilibrium solutions  is an interesting problem. In the case of attractor-bifurcation this problem has already been addressed  by Ma and Wang \cite{MW1} (pp. 155, Theorem 6.1), where one can find an index formula on equilibrium solutions. For the general case treated  here, results in this line will be reported in our forthcoming paper entitled ``Equilibrium index of invariant sets and global static bifurcation for  nonlinear  evolution equations''.
\er

\section{Nontriviality of the Conley Indices of the Bifurcating Invariant Sets}
  Our main goal in this section is to show that the bifurcating  invariant set $K_\lam$ in Theorem \ref{t:2.7} has  nontrivial  Conley index. This result will play a crucial role in establishing our global dynamic bifurcation theorem. However, it may also be of independent interest in its own right.

  Let $m=\mb{dim}\,(X^1)$,  $n=\mb{dim}\,(X^2)$ (\,$n\geq1$\,), and let $K_\lam$ be the  bifurcating  invariant set of $\Phi_\lam$ in Theorem \ref{t:2.7}.
\bt\label{t:4.1}Suppose (H1)-(H4)  are fulfilled (with $\lam_0=0$), and that  $S_0=\{0\}$ is an isolated invariant set of $\Phi_0$.
Then there exists $\ve>0$ such that
\begin{enumerate}
\item[(1)] if $h(\Phi_0,S_0)\ne \Sigma^{m+n}$, then  \be\label{3.28}h(\Phi_\lam, K_\lam)\ne \ol 0,\Hs \lam\in[-\ve,0);\ee
\item[(2)] if $h(\Phi_0,S_0)\ne \Sigma^{m}$, then  \be\label{3.29}h(\Phi_\lam, K_\lam)\ne \ol 0,\Hs \lam\in(0,\ve].\ee
\end{enumerate}

\et

\noindent{\bf Proof.} Let  $U$ be the neighborhood of $0$ given in Theorems \ref{t:2.2} and \ref{t:2.7}.
 Since $S_0$ is an isolated invariant set of $\Phi_0$, we can pick an $\ve>0$ sufficiently small such that $U$ is an isolating neighborhood of the maximal compact invariant set $S_\lam$ of $\Phi_\lam$ for all $\lam\in[-\ve,\ve]$. We may also assume that $U$ and $\ve$ are chosen sufficiently small so that Lemma \ref{l:3.2} applies.
\Vs

(1)\, Assume $h(\Phi_0,S_0)\ne \Sigma^{m+n}.$  Let  $\lam\in[-\ve,0)$. Then by (H1) and (H4),
 $$h(\Phi_\lam,S_0)=\Sigma^{m+n}\ne h(\Phi_0,S_0).$$ and the system bifurcates in $U\sm S_0$ a maximal compact invariant set $K_\lam$.
By  Lemma \ref{l:3.2} one has
$$
h(\Phi_\lam, K_\lam)=h(\Phi_\lam^{12}, K_\lam).
$$
Therefore to prove (\ref{3.28}) we need to check that $h(\Phi_\lam^{12}, K_\lam)\ne \ol 0$.

Choose  an isolating block $N=N_\lam$  of $S_\lam$ in $\cM_\lam^{12}$. Since $S_0$ is a repeller of $\Phi_\lam^{12}$ on $\cM_\lam^{12}$ (by (H4)), one can find an isolating block $N_0$ of $S_0$ in $\cM_\lam^{12}$ (depending upon $\lam$) with $K_\lam\cap N_0=\emp$ such that  $N_0^-=\pa N_0$, where $\pa N_0$ is the boundary of $N_0$ in $\cM_\lam^{12}$.
Then  $M=N\sm\mb{int} N_0$ is an isolating block of $K_\lam$; see Fig.\,5.1.

\vskip-0.8cm
\begin{center}\includegraphics[width=4.3cm]{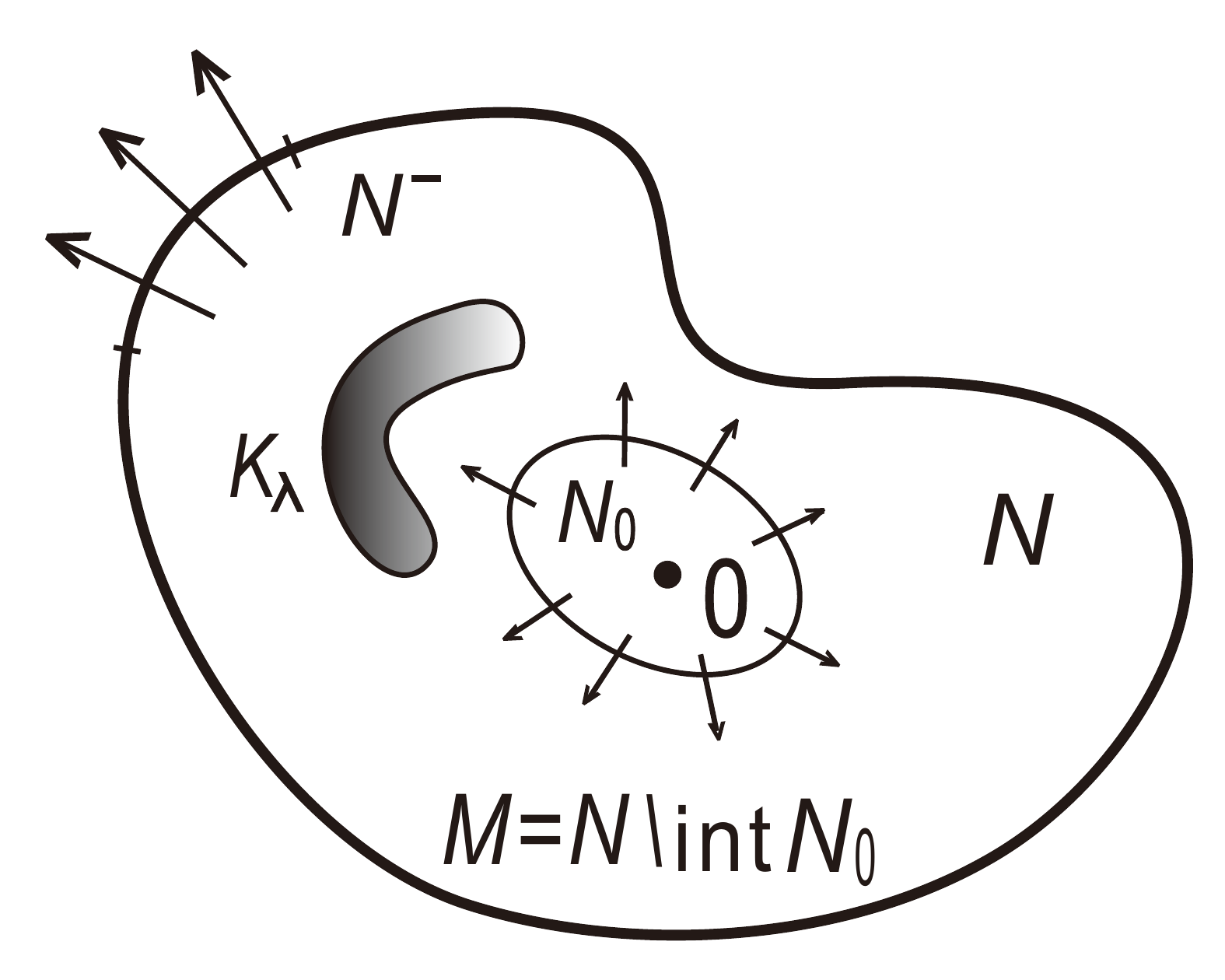} \hspace{1cm} \includegraphics[width=4.3cm]{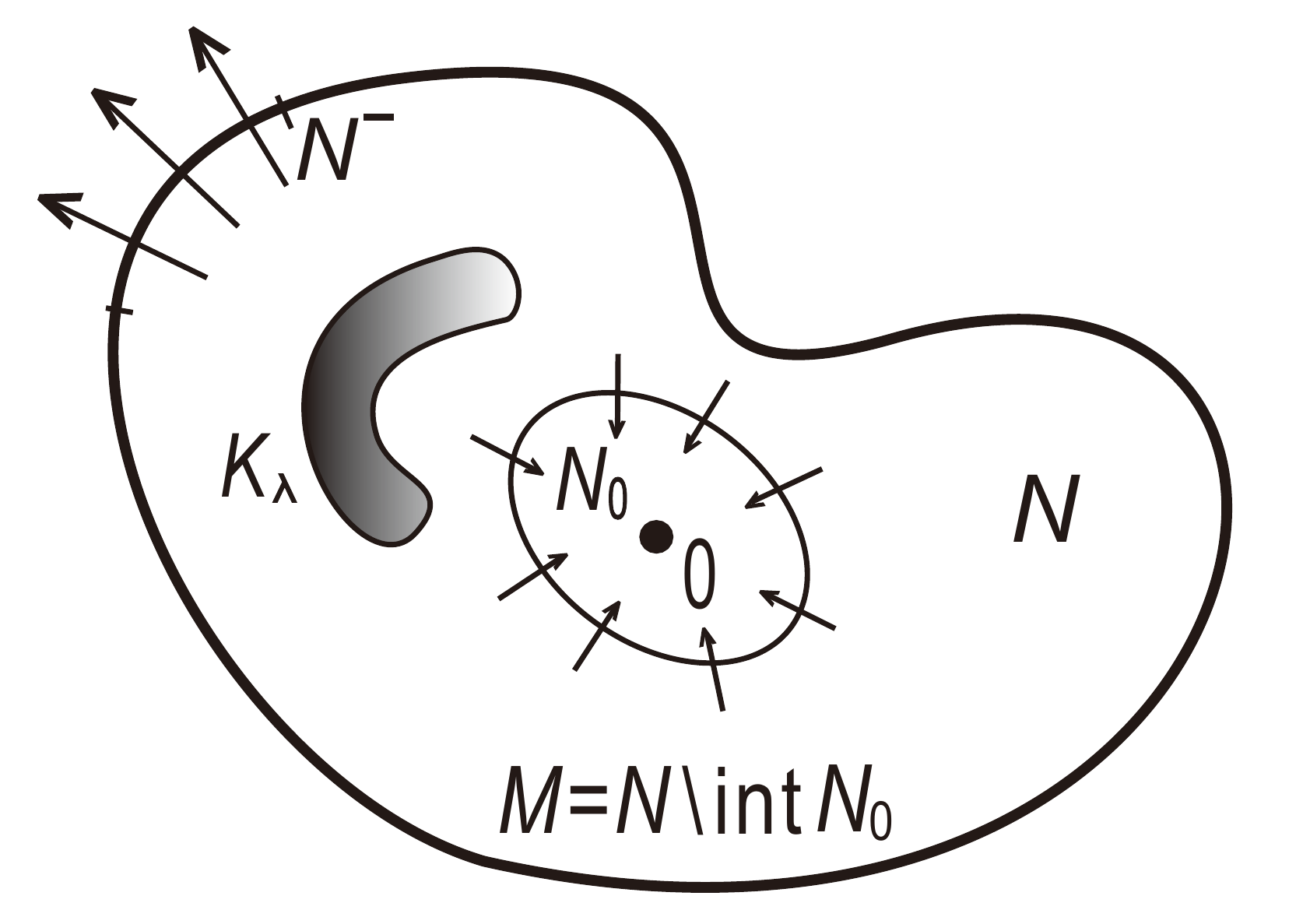}\end{center}
\begin{center}Figure 5.1: \,$\lam<0$ \hspace{3cm} Figure 5.2:  $\lam>0$ \end{center}

As $h(\Phi_\lam,S_0)=\Sigma^{m+n}$, one finds that
\be\label{e5.a1}
S^{m+n}\simeq N_0/\pa N_0= N/M\cong(N/N^-)/\~M,
\ee
where $\~M=\pi_{N^-}(M)$, and $\pi_{N^-}:N\ra N/N^-$ is the projection.
Now let us argue by contradiction and suppose that $h(\Phi_\lam^{12}, K_\lam)=\ol0$. Noticing that $(M/N^-,[N^-])\cong (\~M,[N^-])$ \,(here we have used the same notation $[N^-]$ to denote both the base points in $M/N^-$ and $N/N^-$), we deduce that
$$
[(\~M,[N^-])]=\left[(M/N^-,[N^-])\right]=h(\Phi_\lam^{12}, K_\lam)=\ol0,
$$where ``$\left[\,\.\,\right]$" denotes homotopy type.
This implies that $\~M$ is contractible.

 By a standard argument one can easily show that $\pa N_0$ is a strong deformation retract of $N_0\sm S_0$. Consequently $M$ is a strong deformation retract of $N\sm S_0$. It then follows that $\~M$ is a strong deformation retract of $(N\sm S_0)/N^-$.
Hence by \cite{Ryba} Chap.\,\,I, Pro.\,\,3.6, we deduce  that the pair $(N/N^-,\~M)$  has the homotopy extension property. Further by Lemma \ref{l2.4} and \eqref{e5.a1} it holds that
\be\label{e5.3}
N/N^-\simeq(N/N^-)/\~M\simeq S^{m+n}.
\ee

On the other hand, by the continuation property of the index  we have$$h(\Phi_\lam,S_\lam)= h(\Phi_0,S_0)\ne \Sigma^{m+n},\Hs \lam\in[-\ve,0).$$ Since $h(\Phi_\lam,S_\lam)=h(\Phi_\lam^{12},S_\lam)$, one finds that
  $$h(\Phi_\lam^{12},S_\lam)=\left[(N/N^-,[N^-])\right]\ne \Sigma^{m+n}.$$ This implies that $
N/N^-\not\simeq S^{m+n},
$
which  contradicts (\ref{e5.3}).

\Vs
(2)\, Now consider the case where   $h(\Phi_0,S_0)\ne\Sigma^{m}$.
\vs
Let  $\lam\in(0,\ve]$. Then by (H1) and (H4), we have $h(\Phi_\lam,S_0)=\Sigma^{m}$. Hence
\be\label{e4.0}h(\Phi_\lam,S_0)\ne h(\Phi_0,S_0),\ee
so the system bifurcates in $U\sm S_0$ a maximal compact invariant set $K_\lam\ne\emp$.

Let $S_\lam$ be  the maximal compact invariant set  of $\Phi_\lam$ in $U$. Then $S_\lam\subset \cM_\lam^2$.
By   Lemma \ref{l:3.2} we have
\be\label{e4.10}h(\Phi_\lam,S_\lam)=\Sigma^m\wedge h(\Phi_\lam^2,S_\lam).\ee
On the other hand,  \be\label{4.5}h(\Phi_\lam,S_\lam)= h(\Phi_0,S_0)\ne\Sigma^{m}.\ee
Thus by (\ref{e4.10}) and (\ref{4.5}) one concludes that
\be\label{4.6}
h(\Phi_\lam^2,S_\lam)\ne\Sigma^0.\ee

As  $\Phi_\lam^2$ and the semiflow $\phi_\lam$ generated by  the ODE system (\ref{e:2.9}) on $W$ are conjugate, in the following argument  we {identify} $\Phi_\lam^2$ with $\phi_\lam$, regardless of the conjugacy between them.
By \cite{CE}, Theorem 1.5, one can find  a connected  isolating block $N$ of $S_0$ (with respect to $\phi_0$) with smooth boundary $\pa N$. Further by \cite{CE}, Theorem 1.6, it can be assumed that $\ve$ is sufficiently small so that  $N$ is an  isolating block of $S_{\lam}$ (with respect to $\phi_\lam$) for all $\lam\in(0,\ve]$  with $$B_\lam^-\equiv B^-_0:=N^-,$$ where $B_\lam^-$ denotes the boundary  exit set of $N$ with respect to  $\phi_\lam$.
 We claim that
\be\label{4.7}
N^-\ne\emp.\ee
Indeed, if this was false, $S_0$ would be an attractor of $\phi_0$ in $N$ that attracts $N$. As $S_0$ is a singleton, it follows  that $N$ is contractible. Consequently $$h(\phi_\lam,S_\lam)=h(\phi_0,S_0)=\left[(N/N^-,[N^-])\right]=[(N,\emp)]=\Sigma^0$$
for $\lam\in(0,\ve]$, which  contradicts (\ref{4.6}).

Because $S_0$ is an attractor of $\phi_\lam$ in $W$ for $\lam\in(0,\ve]$ (by (H4)), using appropriate smooth Lyapunov function of $S_0$  one can find an arbitrarily  small isolating block $N_0$ of $S_0$ (depending upon $\lam$) with smooth boundary $\pa N_0$ such that  $N_0^-=\emp$, where $N_0^-$ is the boundary exit set of $N_0$ with respect to $\phi_\lam$.
Note that $M:=N\sm\mb{int} N_0$ is then an isolating block of $K_\lam$ (with respect to $\phi_\lam$) with $M^-=N^-\cup \pa N_0;$ \,see Fig. 5.2.
 We show that 
\be\label{4.8b}CH_*(\phi_\lam, K_\lam)=H_*\(h(\phi_\lam, K_\lam)\)\ne 0,\ee
where $CH_*(\phi_\lam, K_\lam)$ is the homology Conley index of $K_\lam$ with respect to $\phi_\lam$.

First, we infer from  \cite{Ryba}  that the inclusion $M^-\subset M$ has the homotopy extension property. This implies that $M^-$ is a strong deformation retract of one of its neighborhoods in $M$. As $N^-$ and $\pa N_0$ are disjointed compact subsets of $M$,  each of them is a strong deformation retract of a neighborhood of itself in $M$. We collapse $N^-$ and $\pa N_0$ to two distinct  points $z$ and $w$ (see Fig.\,5.3), respectively, and denote $\~M$  the corresponding quotient space. Let $\~M_0=\{z,w\}$. Then
\be\label{4.11}
h(\phi_\lam, K_\lam)= [(M/M^-,[M^-])]=[(\~M/\~M_0,[\~M_0])].\ee

Consider the mapping cone $C_f$ as depicted in Fig.\,5.3, where $f:\~M_0\ra\~M$ is the inclusion. Let $$C\~M_0=(\~M_0\X I)/(\~M_0\X\{1\}).$$ Then $C\~M_0$ is homeomorphic to $I=[0,1]$. Hence one can think of $C_f$ as the space obtained by identifying the end points $0$ and $1$ of $I$ with  $z$ and $w$, respectively, in the  disjoint union of $\~M$ and $I$.
We  observe that $\~M_0$ is a strong deformation retract of an appropriate neighborhood in $\~M$. Consequently  $C\~M_0$ is a strong deformation retract of an appropriate neighborhood in $C_f$. Noticing that $C_f$ is metrizable, by  \cite{Ryba} Chap.\,\,I, Pro.\,\,3.6, we deduce that the inclusion $C\~M_0\subset C_f$  has the homotopy extension property. Since $C\~M_0$ is contractible, by the basic knowledge on homotopy equivalence  (see e.g. \cite{Hat}, Pro.\,\, 0.17), we have
\be\label{4.12}
\~M/\~M_0= C_f/C\~M_0\simeq C_f.
\ee
\vskip-0.8cm
\begin{center}\includegraphics[width=11.5cm]{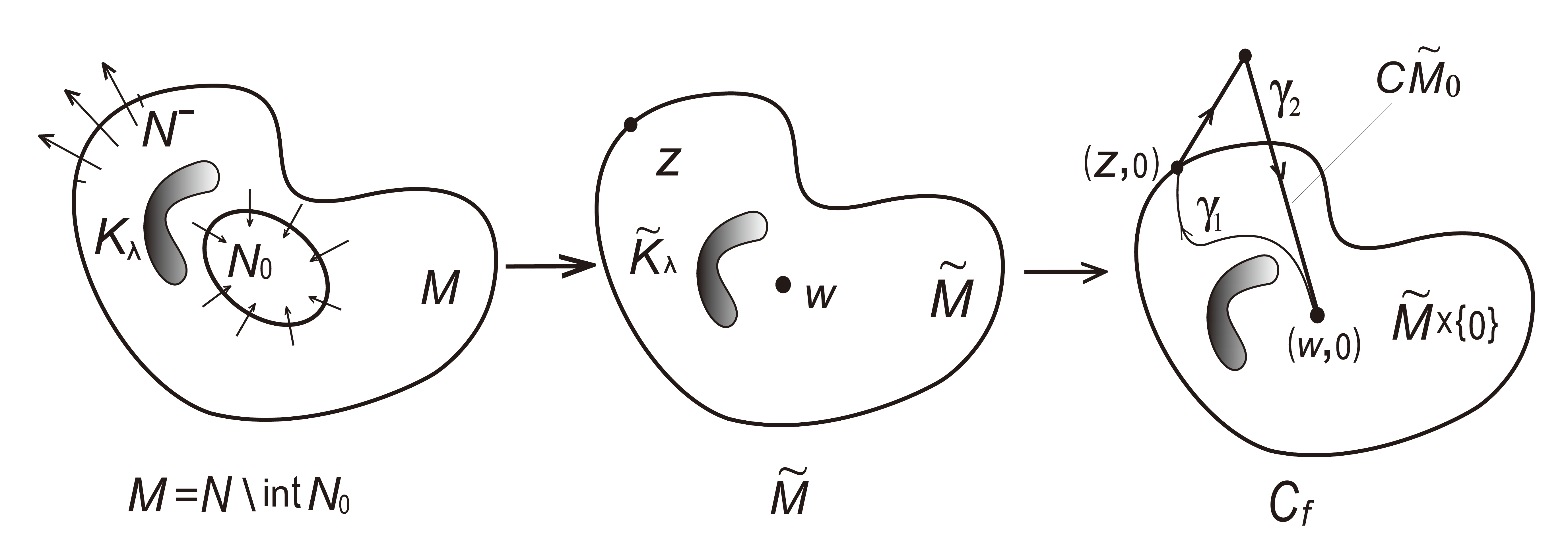}\end{center}
\begin{center}Figure 5.3: \,$M/M^-\simeq\~M/\~M_0\simeq C_f/C\~M_0\simeq C_f$ \end{center}

Because  $M$ is a domain in $E^2$ with smooth boundary, we deduce that $M$ is path-connected. It then follows that $\~M$ is path-connected as well.
Consequently $C_f$ is a path-connected space. Let $\gamma_1$ be a path in $\~M\X\{0\}$ from $(w,0)$ to $(z,0)$ (see Fig.\,5.3), and  $\gamma_2$ be a path in $C_f$ from $(z,0)$ to $(w,0)$ along
$C\~M_0$. Define a closed path $\gamma$ in $C_f$ from $(z,0)$ to $(z,0)$ to be the product $\gamma_1*\gamma_2$ of $\gamma_1$ and $\gamma_2$. Then by a simple continuity argument it can be easily shown that $\gamma$ is not homotopic to  any constant path. Thus the fundamental group $\pi_1(C_f)\ne 0$. Further by some basic knowledge in the theory of algebraic topology we know that $H_1(C_f)\ne 0$. In view of (\ref{4.11}) and (\ref{4.12}) one immediately concludes that $H_1\(h(\phi_\lam, K_\lam)\)\ne 0$. This finishes the proof of (\ref{4.8b}).

\Vs
Now we verify that $h(\Phi_\lam,K_\lam)\ne\ol0$. By Lemma \ref{l:3.2} it suffice to check that $h\(\Phi_\lam^{12},K_\lam\)\ne\ol0$. Suppose the contrary.  Then
we would have $CH_*\(\Phi_\lam^{12},K_\lam\)=0$. Invoking the Poincar$\acute{\mb{e}}$-Lefschetz duality theory on homology Conley index (see McCord \cite{McC}, Theorem 2.1), it then holds that
$CH^*\((\Phi_\lam^{12})^-,K_\lam\)=0$,  where $(\Phi_\lam^{12})^-$ denotes the inverse flow of $\Phi_\lam^{12}$. On the other hand, for $(\Phi_\lam^{12})^-$ we have
$$
h\((\Phi_\lam^{12})^-,K_\lam\)=h\((\Phi_\lam^2)^-,K_\lam\)=h\(\phi_\lam^-,K_\lam\).
$$
(Recall that we identify $\Phi_\lam^2$ with $\phi_\lam$, regardless of the conjugacy between them.)
Hence $CH^*\(\phi_\lam^-,K_\lam\)=0$. Again by the Poincar$\acute{\mb{e}}$-Lefschetz duality theory we find that
$CH_*\(\phi_\lam,K_\lam\)=0$, which contradicts (\ref{4.8b}). \,$\Box$

\section{Global Dynamic Bifurcation}

In this section we establish a  global dynamic bifurcation result.

\subsection{Existence of a local bifurcation branch}
We first prove an existence result for local bifurcation branch.

Set $\cE=E\X\R$, where $E=X^\a$. $\cE$ is equipped with the  metric $\rho$ defined as
$$
\rho\((u,\lam),\,(v,\lam')\)=||u-v||_\a+|\lam-\lam'|, \Hs \A\, (u,\lam),\, (v,\lam')\in\cE.
$$
Let    $\cZ\subset \cE$.  For any $\lam\in\R$,  denote $\cZ_\lam$ the {\em $\lam$-section}  of $\cZ$,
$$
\cZ_\lam=\{u:\,\,(u,\lam)\in \cZ\}.
$$

Let $\~\Phi$ be the {\em skew-product flow} of the family $\Phi_\lam$ ($\lam\in \R$) on $\cE$,
\be\label{e6.1}
\~\Phi(t)(u,\lam)=\(\Phi_\lam(t)u,\,\lam\),\Hs\A\,(u,\lam)\in\cE.
\ee
By the basic theory on abstract evolution equations (see e.g.\,\,\cite{Henry}, Chap. 3 or \cite{Ryba}, Chap. 1, Theorem 4.4 ), one can easily verify  that $\~\Phi$ is {\em asymptotically compact}, i.e., $\~\Phi$ satisfies the hypothesis  (AC) in Section 2.

For each $\lam\in\R$,  denote  $\stackrel{\circ}{\sK}_\lam$ the {\em family of nonempty compact invariant sets $K$ of \,$\Phi_\lam$ with $0\not\in K$}.
Given $\cU\subset \cE$, define
$$\ba{ll}
\sC(\cU)=\ol{\Cup\{K\X\{\lam\}\subset \cU:\,\, K\in \stackrel{\circ}{\sK}_\lam,\,\,\lam\in\R\}}.\ea
$$

\begin{definition}(Bifurcation branch)\, Let $(0,\lam_0)$ be a bifurcation point, and   $\cU\subset \cE$ be a closed neighborhood of  $(0,\lam_0)$.  Then the {bifurcation branch} in $\cU$ from $(0,\lam_0)$, denoted by $\Gamma_\cU(0,\lam_0)$, is defined to be the connected component of $\sC(\cU)$ which contains $(0,\lam_0)$.
\end{definition}


Now we  prove the following interesting  result which ensures  the existence of local bifurcation branch.

\bt\label{t:4.3}Suppose the hypotheses (H1)-(H4) in Theorem \ref{t:2.7} are fulfilled with $\lam_0=0$, and that $S_0=\{0\}$ is an isolated invariant set of $\Phi_0$.
Then there exists $\ve>0$ such that  $$ \Gamma\cap (U\X\{\pm\ve\})\ne\emp,$$ where $\Gamma=\Gamma_\cU(0,0)$, and $\cU=U\X[-\ve,\ve]$.

\et

\Vs\noindent{\bf Proof.} Let $U$  be the neighborhood of $0$ given in Theorem \ref{t:2.7}, and let $S_\lam$ be the maximal compact invariant set of $\Phi_\lam$ in $U$.
Choose an  $\ve>0$ such  that the assertions in Theorem \ref{t:4.1} 
hold. Let $K_\lam$ be the maximal compact invariant set of $\Phi_\lam$ in $U\sm S_0$.  Since $\lim_{\lam\ra0}d_{\mbox{\tiny H}}(K_\lam,0)=0$, we may also assume $\ve$ is sufficiently small so that there exists $r>0$ such that \be\label{e:5.3}\mB(K_\lam,r)\subset U,\Hs\A\,\lam\in[-\ve,\ve].\ee
We show that $\ve$ fulfills the requirement of the theorem.

For definiteness, by Theorem \ref{t:4.1} it can be  assumed that
\be\label{e6.2}
h(\Phi_\lam, K_\lam)\ne \ol 0
\ee
 for $\lam\in(0,\ve]$. We check that  $$ \Gamma\cap (U\X\{\ve\})\ne\emp,$$ thus completing  the proof of the theorem.

\vs We  first prove that for any $0<\mu<\ve$, $\sC(\cU_\mu)$ has a connected component $\cZ$ such that
\be\label{e:4.10}\ba{ll}
\cZ\cap \(U\X\{\mu\}\)\ne\emp\ne \cZ\cap \(U\X\{\ve\}\),\ea
\ee
where $\cU_\mu=U\X[\mu,\ve]$. For this purpose, let us  first verify that
$$\ba{ll}\sC(\cU_\mu)={\Cup_{\mu\leq\lam\leq\ve}\,K_\lam\X\{\lam\}}:=\cK.\ea$$
 Indeed, we infer from the maximality of $K_\lam$  in $U\sm S_0$ that  $\sC(\cU_\mu)=\ol{\cK}.$ On the other hand, it is clear that $\cK$ is invariant under the skew-product flow $\~\Phi$. Hence by  asymptotic  compactness of $\~\Phi$ we deduce that $\cK$ is pre-compact. Further  by upper semicontinuity of $K_\lam$ in $\lam$ one can easily verify  that $\cK$ is  closed. Thus $\cK$ is compact. Consequently $\sC(\cU_\mu)=\ol\cK=\cK.$

The compactness of $\cK$ also implies
\be\label{e:5.5}
d(0,K_\lam)\geq 2\eta,\Hs\A\,\lam\in[\mu,\ve],\ee
 where  $\eta>0$ is a positive number independent of $\lam$.

In what follows we argue by contradiction and suppose that (\ref{e:4.10}) fails to be true.
Then for any connected component $\cZ$ of $\sC(\cU_\mu)$  one  has
$$\ba{ll}\mb{either}\,\,\,\cZ \cap \(U\X\{\mu\}\)=\emp,\hs \mb{or }\,\cZ \cap\( U\X\{\ve\}\)=\emp.\ea$$
If there are only a finite number of components, then each component $\cZ$ is isolated in $\cU$. Because the $\lam$-section $\cZ_\lam$ of  $\cZ$ is empty when $\lam$ is close to either $\mu$ or $\ve$, by the continuation property of Conley index we see that $h(\Phi_\lam,\cZ_\lam)\equiv \ol0.$ Consequently the ``sum'' of these indices equals $\ol0$. This contradicts \eqref{e6.2} and justifies  (\ref{e:4.10}), as the union of ${\cZ_\lam}'s$ is precisely $K_\lam$. However, in general there is also the possibility that $\sC(\cU_\mu)$ may contain  infinitely many components. We will employ the Separation Lemma given in Section 2 to overcome this difficulty.

Set $\cO_\mu=\cU_\mu\sm\(\mB(0,\eta)\X[\mu,\ve]\)$. Then clearly $\sC(\cO_\mu)=\sC(\cU_\mu)$. Denote $\sF$  the family of connected components of $\sC(\cO_\mu)$. By (\ref{e:5.3}) and (\ref{e:5.5}) we see  that $\cO_\mu$ is a neighborhood of $\cZ$ in the space  $$\cH=E\X[\mu,\ve]$$ for each $\cZ\in \sF$.  This allows us to pick for each $\cZ\in \sF$ a closed neighborhood $\Omega_\cZ$ in $\cH$ with $\Omega_\cZ\subset \cO_\mu$ such that
if $\cZ\cap \(U\X\{\sig\}\)=\emp$ (where  $\sig=\mu$ or $\ve$), then \be\label{Om}\ba{ll}\Omega_\cZ\cap \(U\X\{\sig\}\)=\emp;\ea\ee see Fig.\,\,6.1.
\vskip-0.3cm\begin{figure}[h]\centering\includegraphics[width=6cm]{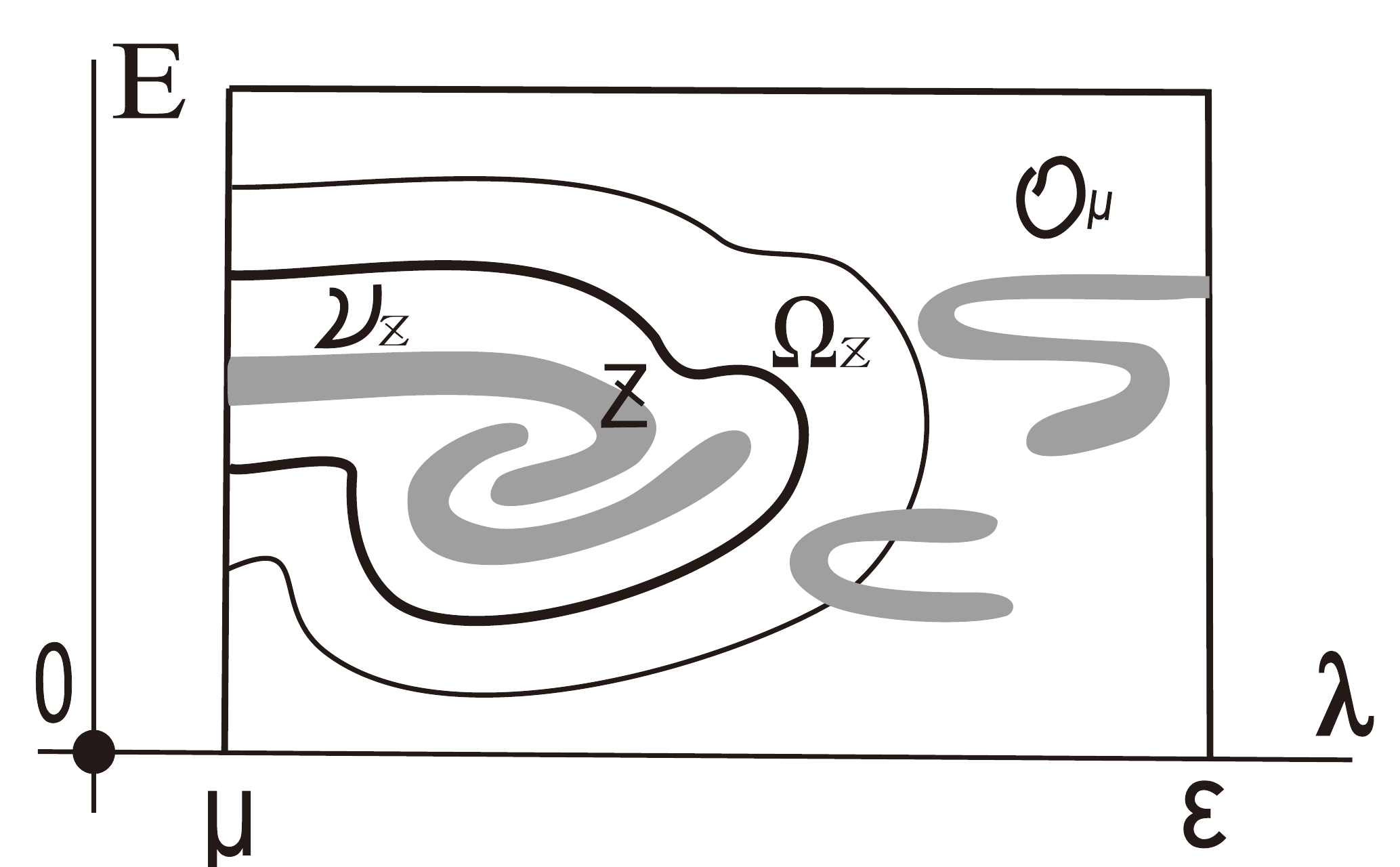}
\caption{Separating neighborhoods  of $\cZ$ in $\cH$}\label{fg4-1}\end{figure}

For any $\cO\subset \cH$,   denote  $\pa_\cH\cO$ the boundary of $\cO$ in  $\cH$.
Given $\cZ\in\sF$, set
$$\ba{ll}
\mathfrak{B}=\Cup\{\cF\in \sF:\,\,\cF\cap \,\pa_\cH\Omega_\cZ\ne\emp\}, \hs \mathfrak{D}=\Cup\{\cF\in\sF:\,\,\cF\cap\, \Omega_\cZ\ne\emp\},\ea
$$
 We claim that both $\mathfrak{B}$ and $\mathfrak{D}$ are closed. Indeed, if $b\in \ol{\mathfrak{B}}$, then there exists a sequence $b_k\in \mathfrak{B}$ such that $b_k\ra b$.
We may assume that $b_k\in \cF_k$ for some $\cF_k\in \sF$ with $\cF_k\cap \pa_\cH\Omega_\cZ\ne\emp$. By Lemma \ref{l:2.3} we deduce that
there exists a subsequence of $\cF_k$, still denoted by $\cF_k$, such that
$$
\lim_{k\ra\8}\de_{\mbox{\tiny H}}(\cF_k,\cF_0)=0.
$$
One trivially checks that $\cF_0$ is connected and contained in $\sC(\cO_\mu)$; moreover, $\cF_0\cap \pa_\cH\Omega_\cZ\ne\emp$. Since
$b\in\cF_0$, we conclude that $b\in \mathfrak{B}$. Hence $\mathfrak{B}$ is closed. Likewise it can be shown that $\mathfrak{D}$ is closed.

Note that  $\cZ\cap\mathfrak{B}=\emp$. Since  $\cZ$ does not intersect any other connected component of $\mathfrak{D}$, by Lemma \ref{l:2.2} there exist two disjoint closed subsets $\cK_1$ and $\cK_2$ of $\mathfrak{D}$ such that $\mathfrak{D}=\cK_1\cup\cK_2$, and $$\cZ\subset \cK_1,\hs \mathfrak{B}\subset \cK_2\,.$$ It is clear that $\cK_1$ is contained in the interior   of $\Omega_\cZ$ relative to $\cH$.

Take a positive number $\de_\cZ$ with  $$\de_\cZ<\frac{1}{8}\min\(d(\cK_1,\cK_2),\,d(\cK_1,\pa_\cH\Omega_\cZ)\).$$
Let $\cV_\cZ=\mB_{\cH}(\cK_1,4\de_\cZ)$ be the $4\de_\cZ$-neighborhood of $\cK_1$ in $\cH$.  Then  $\cV_\cZ\subset \Omega_\cZ$, and
\be\label{e:5.6}\ba{ll}\mB_{\cH}(\pa_\cH\cV_\cZ,2\de_\cZ)\cap \sC(\cO_\mu)=\emp.\ea\ee
By the compactness of $\sC(\cO_\mu)$ there exist a finite number of $\cZ\in\sF$, say, ${\cZ_1},\cdots,\cZ_l$,  such that
$\sC(\cO_\mu)\subset\Cup_{1\leq k\leq l}\cV_{\cZ_k}\,.$
Set
$$\ba{ll}
\cW_k=\cV_{\cZ_k}\sm\(\ol\cV_{\cZ_1}\cup\cdots\cup  \ol\cV_{\cZ_{k-1}}\),\Hs k=1,2,\cdots,l.\ea
$$
Then ${\cW_k}'s$ are disjoint open subsets of $\cH$. One can easily check that
 \be\label{e6.7}\ba{ll}\pa_\cH\cW_k\subset \Cup_{1\leq i\leq k}\pa_\cH\cV_{\cZ_i}.\ea\ee
Thus we deduce that $\sC(\cO_\mu)\subset\Cup_{1\leq k\leq l}\cW_k\,.$
\vs
Let $\cS_k=\sC(\cO_\mu)\cap\cW_k$.  We claim that
\be\label{e:5.7}
d\(\cS_k,\pa_\cH\cW_k\)>0.
\ee
Indeed,
 if $w\in \cS_k$ then by (\ref{e:5.6}) we have
$$
d(w,\pa_\cH\cV_{\cZ_i})\geq 2\de_{\cZ_i}\geq 2\min_{1\leq j\leq l}\de_{\cZ_j}:=\de_0>0,\Hs 1\leq i\leq l,$$
and  the conclusion follows from (\ref{e6.7}).
\vs
 It follows by (\ref{e:5.7}) that $\cS_k=\sC(\cO_\mu)\cap\ol\cW_k$. Hence $\cS_k$ is compact. It can be easily seen that $\cS_k$ is the maximal compact invariant set
 of $\~\Phi$ in $\ol\cW_k$.
 Since $\ol\cW_k$ is a neighborhood of $\cS_k$ in $\cH$, by Theorem \ref{t:2.14} we have
 \be\label{e:4.12}
h(\Phi_\lam,\cS_{k,\lam})\equiv \mb{const.}, \Hs \lam\in [\mu,\ve],
\ee where $\cS_{k,\lam}$ is the $\lam$-section of $\cS_k$.
On the other hand, by (\ref{Om}) we have either $\cS_{k,\mu}=\emp$, or $\cS_{k,\ve}=\emp$. Hence by (\ref{e:4.12}) it holds  that
\be\label{e:4.12b}
h(\Phi_\lam,\cS_{k,\lam})\equiv \ol0, \Hs \lam\in [\mu,\ve],
\ee

Now by  (\ref{e:4.12b}) we conclude  that
$$
h(\Phi_\lam,K_{\lam})=h(\Phi_\lam,\cS_{1,\lam})\vee \cdots \vee h(\Phi_\lam,\cS_{l,\lam})
=\ol0.
$$
This contradicts (\ref{e6.2}) and completes the proof of (\ref{e:4.10}).
\vs
We are now ready to complete the proof of the theorem. Take a sequence of positive numbers  $\mu_k\ra 0$. For each $\mu_k$,  pick  a connected component $\cZ_k$ of $\sC(\cO_{\mu_k})$ such that
$$
\ba{ll}\cZ_k\cap \(U\X\{\mu_k\}\)\ne\emp\ne \cZ_k\cap \(U\X\{\ve\}\).\ea
$$
By Lemma \ref{l:2.3} we may assume  that
$$
\lim_{k\ra\8}\de_{\mbox{\tiny H}}(\cZ_k,\cZ_0)=0.
$$
Then $\cZ_0$ is a continuum in $\sC(\cU)$ with $(0,0)\in \cZ_0$ and $\cZ_0\cap (U\X\{\ve\})\ne\emp$.\, $\Box$

\subsection{Global bifurcation}

For the sake of convenience in statement, we make a convection that $\8\in \pa \Omega$ if $\Omega$ is an unbounded subset of $\cE$.

The main result in this section is  the following
theorem.
\bt\label{gbt}(Global dynamic bifurcation) Assume that the hypotheses in Theorem \ref{t:2.7} are fulfilled.
Let $\Omega\subset \cE$ be a closed neighborhood of the bifurcation point $(0,0)$. Suppose that $S_0=\{0\}$ is an
isolated invariant set of  $\Phi_0$.

Let $\Gamma=\Gamma_\Omega(0,0)$.
Then one of the following cases  occurs.
\begin{enumerate}
\item[(1)] $\Gamma\Cap \pa\Omega\ne\emp$; see Fig.  6.3.
\vs
\item[(2)] $0\in\ol{ \Gamma_0\sm\{0\}}$, where $\Gamma_0$ is the $0$-section of $\Gamma$; see Fig.  6.4.
\vs
\item[(3)] There exists
$\lam_1\ne0$  such that $(0,\lam_1)\in\Gamma$\,; see Fig.  6.5. \end{enumerate} \et


\begin{center}\includegraphics[width=4.3cm]{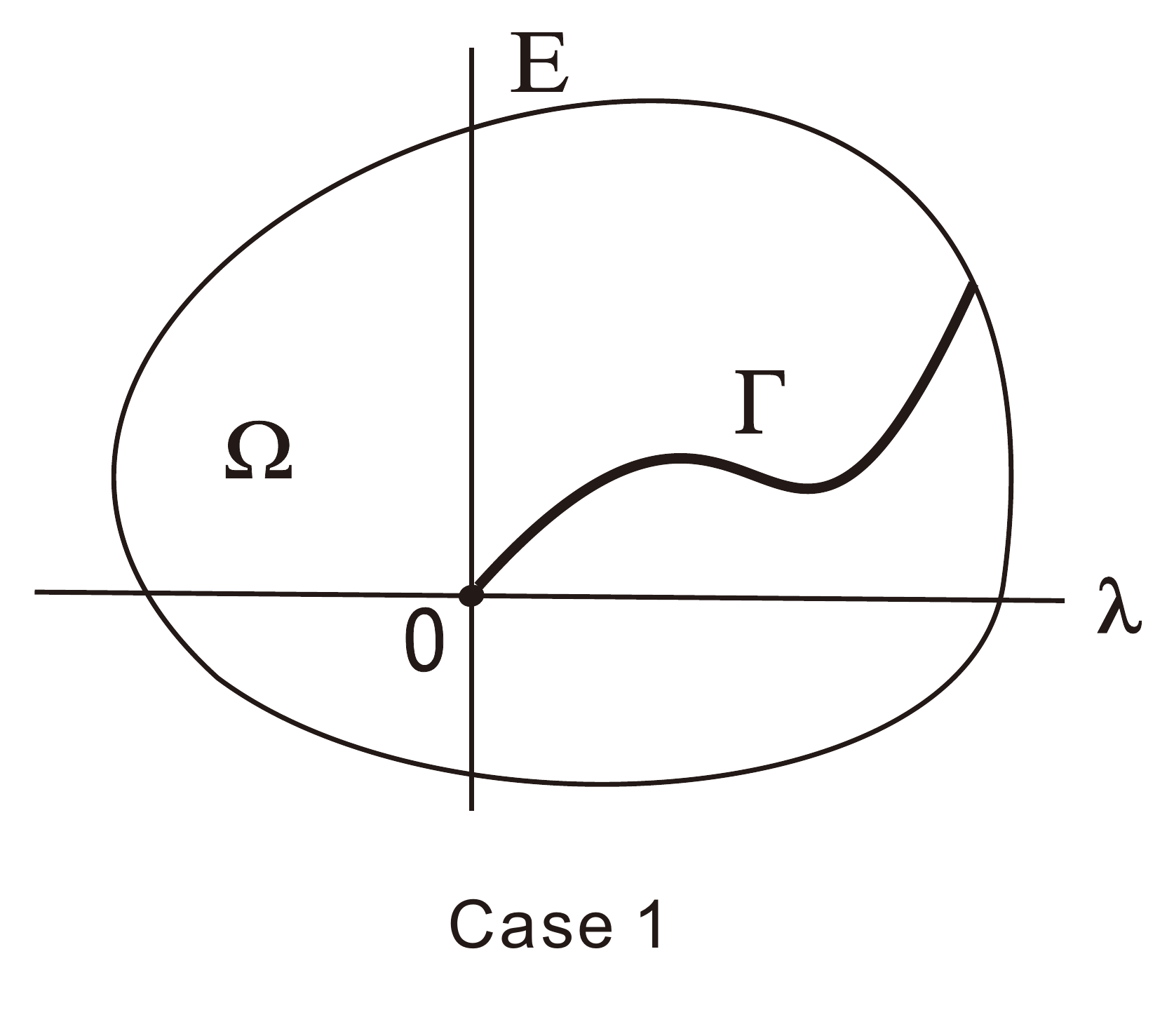} \hspace{1cm} \includegraphics[width=4.3cm]{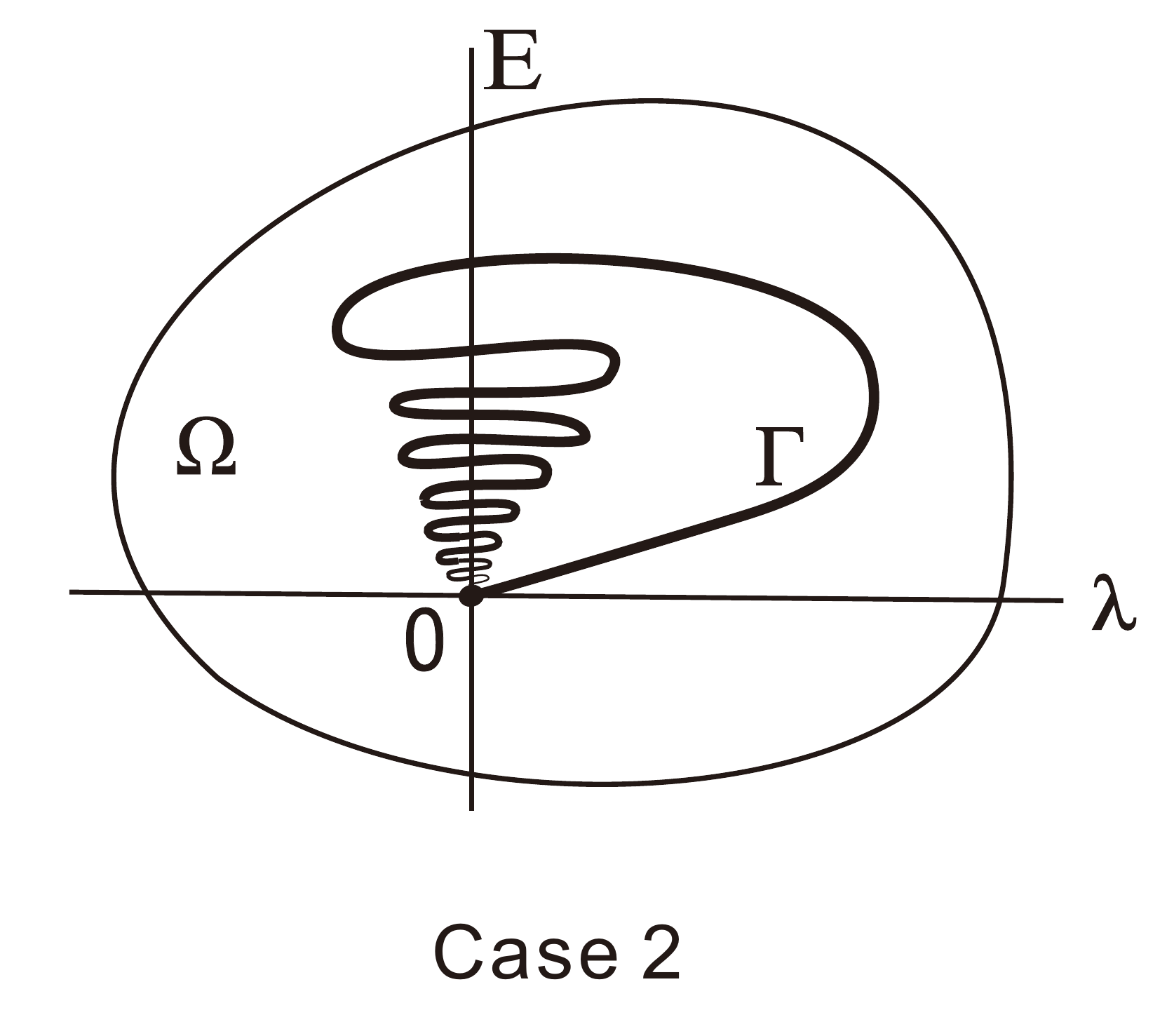}\end{center}
\begin{center}Figure 6.3: \,$\Gamma\Cap \pa\Omega\ne\emp$ \hspace{2.5cm} Figure 6.4: \,$0\in\ol{ \Gamma_0\sm\{0\}}$ \end{center}
\vs
\begin{center}\includegraphics[width=4.3cm]{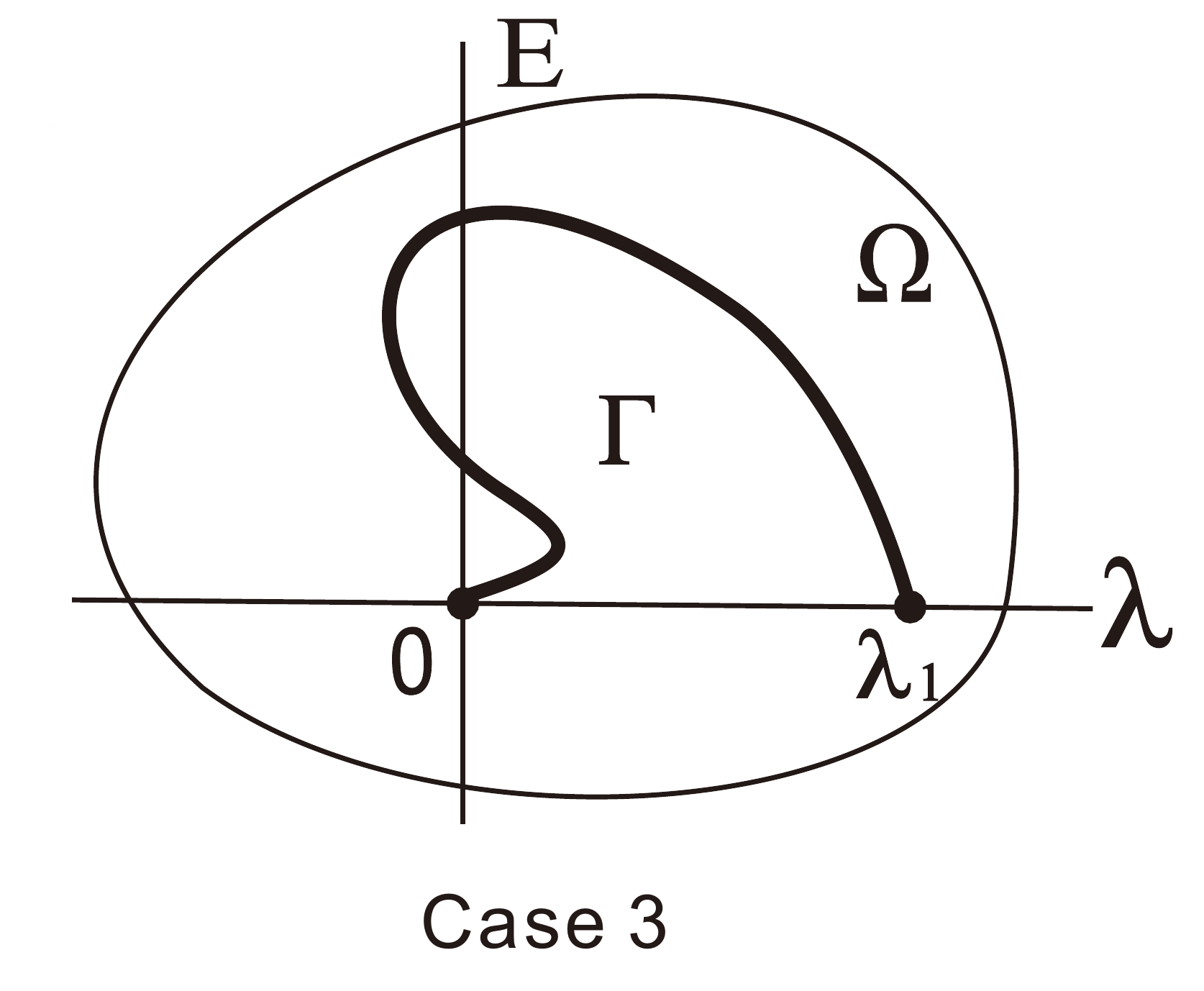} \hspace{1.2cm} \includegraphics[width=4.3cm]{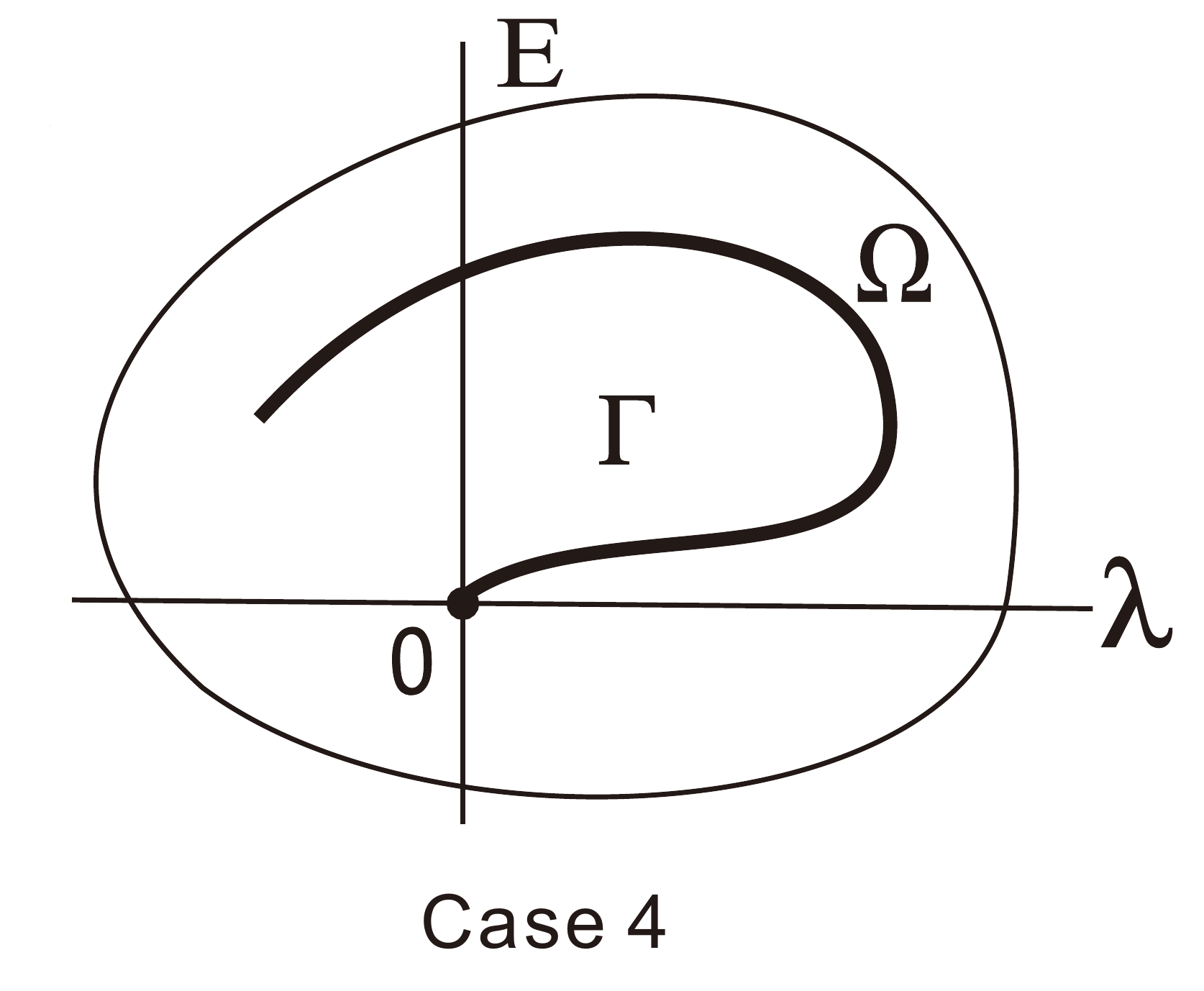}\end{center}
\begin{center} \Hs Figure 6.5: \,$(0,\lam_1)\in\Gamma$ \hspace{2cm} Figure 6.6: This case never occurs\end{center}

\noindent{\bf Proof.} We argue by contradiction and suppose that none of the cases  (1)-(3) occurs. Then $\Gamma$ is  a
bounded closed subset of $\cE$  contained in the
interior  of $\Omega$ as depicted in  Fig.\,\,6.6. It is easy to see  that  $\Gamma$ is invariant under the skew-product flow $\~\Phi$. Hence by asymptotic compactness of $\~\Phi$ we deduce that $\Gamma$ is compact.

Since $0\not\in\ol{ \Gamma_0\sm S_0}$, we can write $\Gamma_0$ as $\Gamma_0=S_0\cup A_0$, where $A_0$ is a compact  invariant set of  $\Phi_0$ with $A_0\cap S_0=\emp$. We only consider the case where $A_0\ne\emp$. The argument for the case where $A_0=\emp$ is a slight modification of that of the former one.

Let $U\subset E$ and $\ve>0$ be  as  in Theorem \ref{t:4.3}.
Then 
the system $\Phi_\lam$ bifurcates, say, for each $0<\lam\leq \ve$, a nonempty maximal compact invariant set $K_\lam$ in $U\sm S_0$ with
\be\label{e5.11}\lim_{\lam\ra0}d_{\mbox{\tiny H}}(K_\lam,S_0)=0\ee and
\be\label{e:5.10}h(\Phi_\lam,K_\lam)\ne\ol0,\Hs \A\,\lam\in (0,\ve].\ee
Pick a closed neighborhood $V$ of $0$ with $V\subset U$ and $$d\(A_0,\,V\):=\sig_0>0.$$
By (\ref{e5.11}) we can further restrict $\ve$ sufficiently small so that for some   $r_0>0$,
$$
\mB(K_\lam,r_0)\subset V,\Hs \A\, \lam\in(0,\ve].
$$

By the compactness of $\Gamma$  it is easy to verify that the $\lam$-section  $\Gamma_\lam$ of $\Gamma$ is   upper semicontinuous in $\lam$.
Let $$\cZ=\Gamma\cap(V\X[0,\ve]).$$ Then $d_{\mbox{\tiny H}}(\cZ_\lam,S_0)\ra0$ as $\lam\ra 0$. As $A_0\cap S_0=\emp$,
it also holds  that $$\lim_{\lam\ra0}d_{\mbox{\tiny H}}(A_\lam,A_0)\ra 0,$$ where $A_\lam=\Gamma_\lam\sm\cZ_\lam$ ($\lam\in[0,\ve]$).
Thus there exist $\eta_0>0$ and $0<\ve'\leq\ve$ such that \be\label{e:5.9}\ba{ll}\ol\mB(\cZ_\lam,\eta_0)\subset V,\hs
\ol\mB(A_\lam,\eta_0)\cap V=\emp\ea
\ee
for all $\lam\in[0,\ve']$.
Note that  both $\cZ_\lam$ and $A_\lam$ are compact invariant sets of $\Phi_\lam$.


\vs

Let
$M_0=\Cup_{\lam\geq\ve'}\Gamma_\lam$. It can be easily shown that  $M_0$ is a compact subset of $E$. Clearly $0\not\in M_0$, hence
\be\label{e6.8}d(0,M_0):=\de_0>0.\ee
Fix a number $0<r<\frac{1}{3}\min\(\eta_0,\de_0\)$.
Utilizing the Separation Lemma, by a similar argument as in
the proof of Theorem \ref{t:4.3} we can find a closed
neighborhood $\cO$ of $\Gamma$ with $\cO\subset \mB_\cE(\Gamma,r)$
such that
\be\label{5.10}\ba{ll}
\sC(\Omega)\cap \pa\cO=\emp.\ea
\ee
Here $\mB_\cE(\Gamma,r)$ denotes the $r$-neighborhood of $\Gamma$ in $\cE$.  By the choice of $r$ it can be easily seen that if $\lam\in(0,\ve']$ then
$$\cO_\lam\subset  \ol\mB(\cZ_\lam,\eta_0)\cup\ol\mB(A_\lam,\eta_0);$$
see Fig.\,6.7. Set
$$\ba{ll}
G_\lam=\cO_\lam\cap \ol\mB(\cZ_\lam,\eta_0),\hs H_\lam=\cO_\lam\cap
\ol\mB(A_\lam,\eta_0).\ea$$ By  (\ref{e:5.9}) we have
\be\label{e:4.17}\ba{ll}\cO_\lam=G_\lam\cup H_\lam,\hs G_\lam\cap
H_\lam=\emp\ea\ee for $\lam\in(0,\ve']$.

We claim that there exists  $\sig>0$ such that $$\mB_{\sig}\subset G_\lam$$ for all $\lam$ sufficiently small,  where (and below) $\mB_R$ denotes the ball in $E$ centered at $0$ with radius $R$.
Suppose the contrary. There would exist sequences $\lam_k\ra0$ and $x_k\in \pa G_{\lam_k}$  such that $x_k\ra 0$. Noticing that
$(x_k,\lam_k)\in\pa\cO$, one concludes  that $(0,0)\in \pa\cO$, a
contradiction!

By (\ref{e5.11}) one can find  a  number
$0<\mu\leq{\ve'}/{2}$ such that \be\label{e:4.15} K_\lam\subset
\mB_\sig\subset G_\lam\,,\Hs\,\A \lam\in(0, 2\mu].\ee  Using the  upper semicontinuity of $K_\lam$ in $\lam$  (see Theorems {\ref{t:2.7}}) one can easily show that $F=\Cup_{\mu\leq\lam\leq\ve'}K_\lam$ is  closed in $E$. Because   $\cF=\Cup_{\mu\leq\lam\leq\ve'}K_\lam\X\{\lam\}$ is invariant under the system $\~\Phi$, by asymptotic compactness of $\~\Phi$  we deduce  that $\cF$ is pre-compact in $\cE$. It then follows that $F$ if compact in $E$. Hence
$$
d(0,F):=d_0>0.
$$

Take a $\Lam>0$  such that $\cO\subset
E\X(-\Lam,\Lam)$. Let $\rho$ be a positive number with
$\rho<\rho_0:=\frac{1}{2}\min(d_0,\de_0)$, where $\de_0$ is the number given in \eqref{e6.8}. Set
$$\ba{ll}
\cV=\cO\cap \cH,\hs
\cW=\cV\sm\(\mB_\rho\X[\mu,\Lam]\),\ea
$$where $\cH=E\X[\mu,\Lam]$; see Fig.\,\,6.8. Clearly $\cV$ is closed in $\cH$. Since $\mB_\rho\X[\mu,\Lam]$ is open in $\cH$, we see that $\cW$ is closed in $\cH$ as well.
We claim that \be\label{5.16}\sC(\cW)=\sC(\cV):=\sC,\Hs \A\,\rho<\rho_0.\ee To see this, by definition it suffices to show that if $\lam\in [\mu,\Lam]$, then  any compact invariant set $M$ of $\Phi_\lam$ in
$\cV_\lam\sm S_0$ is necessarily  contained in $\cW_\lam$.
\vskip-0.8cm
\begin{center}\Hs\includegraphics[width=6cm]{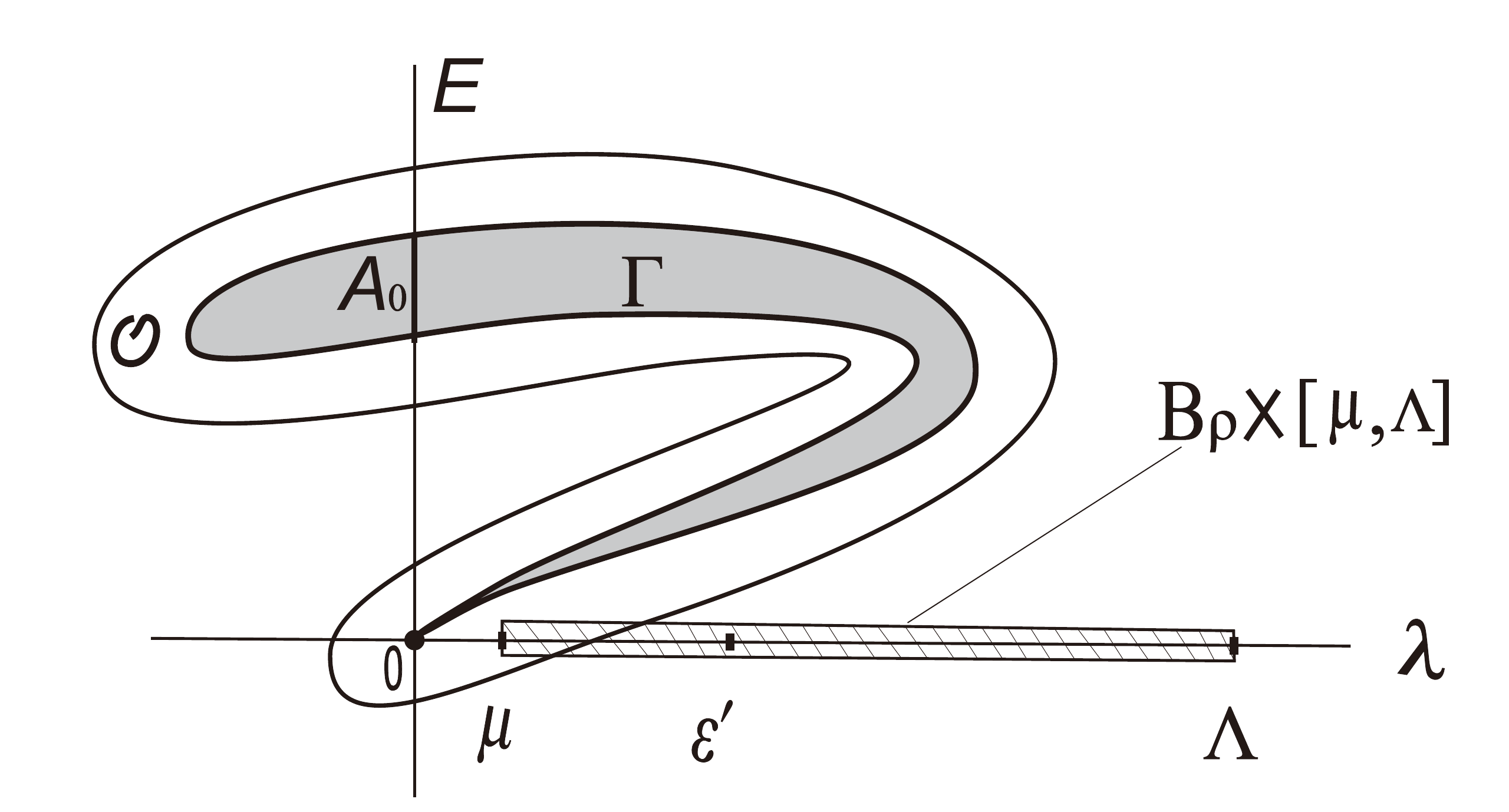}\includegraphics[width=6cm]{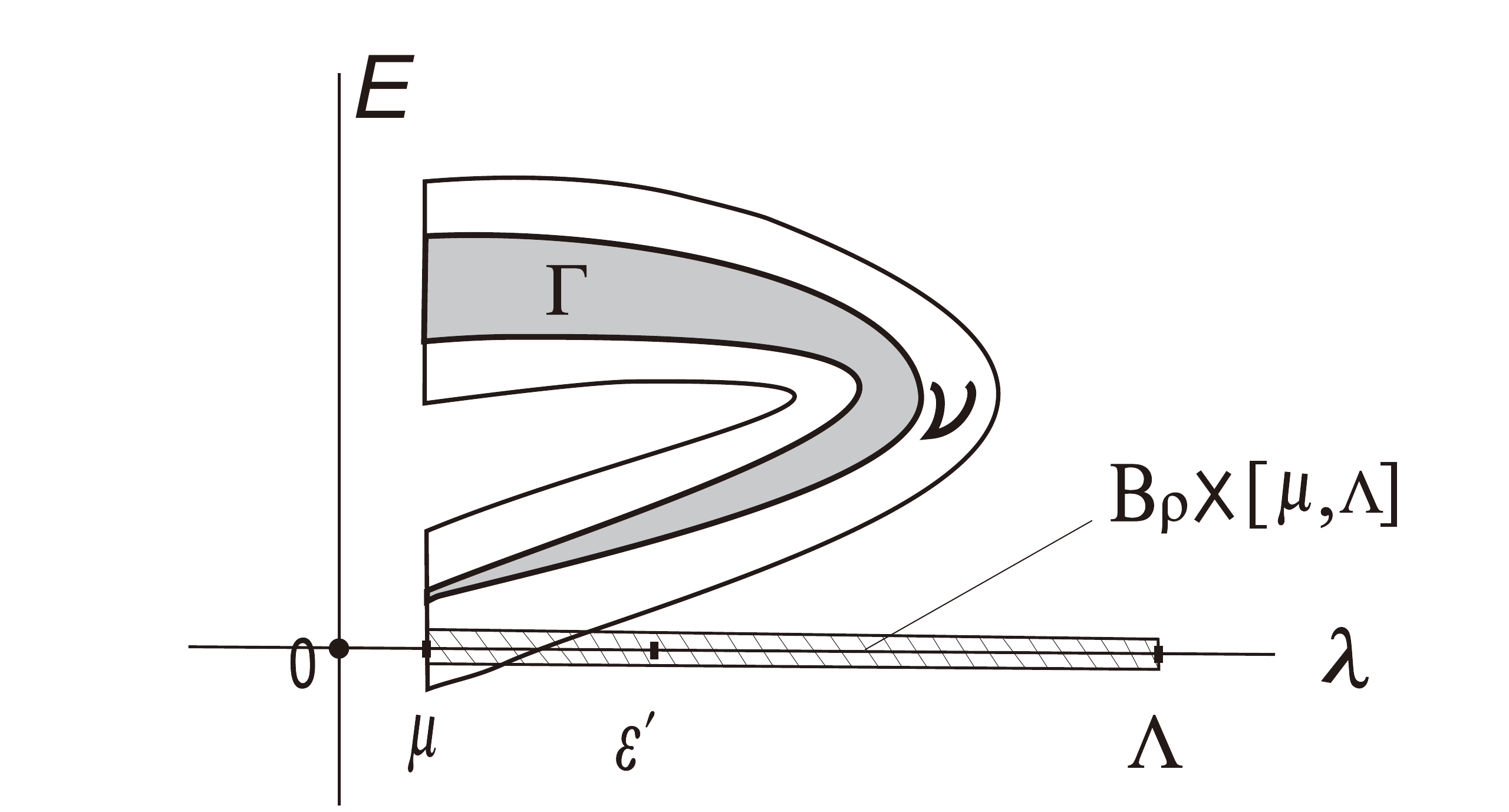}\\[2ex] Figure 6.7 \hspace{5cm} Figure 6.8 \end{center}

We first consider  $\lam\in[\mu,\ve']$. By
(\ref{e:5.9}) and the choice of $r$, we have   $$\ba{ll}M\subset \cV_\lam\subset \ol\mB(A_\lam,\eta_0)\cup V.\ea$$
Clearly  $M\cap\ol\mB(A_\lam,\eta_0)\subset\cW_\lam$. We observe  that  $M\cap V$ is a compact invariant set of $\Phi_\lam$ in $V\sm S_0$.
Therefore by the maximality of $K_\lam$ in $V\sm S_0$ one has  $M\cap V\subset K_\lam$.
Because $K_\lam\cap \mB_{2\rho}=\emp$ for $\mu\leq \lam\leq \ve'$,  by the definition of $\cW$ we see that $M\cap V\subset
\cW_\lam$. Thus $M\subset \cW_\lam$.

 Now assume that  $\lam>\ve'$. Then  by the choice of
$\rho$ we find that $\cO_\lam\cap \mB_\rho=\emp$; see Fig.\,\,6.8. It follows that
$\cV_\lam=\cW_\lam$.  This finishes the proof of what we desired. Hence (\ref{5.16}) holds true.

We show that $\cV$ is a neighborhood of $\sC$ in $\cH:=E\X[\mu,\Lam]$. Suppose the contrary. Then $\sC\cap \pa_{\cH}\cV\ne\emp$, where $\pa_\cH\cV$ denotes the boundary of $\cV$ relative to $\cH$. Noticing that $$\ba{ll}\pa_{\cH}\cV=\pa_{\cH}(\cO\cap \cH)\subset\pa\cO\cap \cH,\ea$$ we have $$\ba{ll}\sC\cap \pa \cO=\sC\cap \(\pa\cO\cap \cH\)\supset \sC\cap \pa_{\cH}\cV \ne\emp.\ea$$ This contradicts  (\ref{5.10}).

By (\ref{5.16})  we can  fix a $\rho>0$ sufficiently small so that $\cW=\cV\sm\(\mB_\rho\X[\mu,\Lam]\)$ is a neighborhood of $\sC$ in $\cH$. By the definitions of $\sC=\sC(\cW)$ and the skew-product flow one can easily see that $\sC$ is  the maximal compact invariant set of  $\~\Phi$ in $\cW$. Hence $\cW$ is an isolating neighborhood of $\sC$ in $\cH$.
  It then follows by Theorem \ref{t:2.14} that
 \be\label{e:4.16} h(\Phi_\lam,\sC_\lam)\equiv  h(\Phi_{\Lam},\sC_{\Lam})=  h(\Phi_{\Lam},\emp)=\ol0,\Hs \lam\in[\mu, \Lam].\ee

On the other hand, if $\mu\leq\lam\leq2\mu$ then by (\ref{e:4.15}) and the choice of $\rho$, we find that  $\tilde{G}_\lam:=G_\lam\sm\mB_\rho$ is a
neighborhood of $K_\lam$. Since $K_\lam$ is the maximal compact invariant of $\Phi_\lam$ in $V\sm S_0$ (and hence in $\tilde{G}_\lam$), we infer from
(\ref{e:4.17})   that   $\sC_\lam\sm K_\lam$ is necessarily contained in $ H_\lam$ (note that $\cW_\lam=\~G_\lam\cup H_\lam$).  Thus
$$
h(\Phi_\lam,\sC_\lam)=h(\Phi_\lam,K_\lam)\vee
h(\Phi_\lam,\sC_\lam\sm K_\lam).
$$
 (\ref{e:4.16}) then implies that   $$h(\Phi_\lam,K_\lam)=\ol0.$$  This contradicts  (\ref{e:5.10}), which completes the proof of the theorem. \,$\Box$

\section{An Example}
In this  section we give an example to illustrate our theoretical results by considering the well-known Cahn-Hilliard equation describing the  spinodal decomposition.

 The  nondimensional form of the equation reads (see \cite{MW2})
\be\label{e:6.2}\left\{\ba{lll}
u_t+\Delta^2 u+\lam\De u=\De(b_2 u^2+b_3 u^3),\Hs &(x,t)\in \Omega\times R^+,\\[1ex]
\frac{\pa u}{\pa \nu}=\frac{\pa (\De u)}{\pa \nu}=0,\Hs &(x,t)\in \partial\Omega\times R^+,\\[1ex]
m(u)=0,\ea\right.
\ee
where   $\Omega\subset R^d$ ($d\leq 3$) is a bounded domain with  smooth boundary
$\partial\Omega$,  $b_3>0$, and $$
m(u)=\frac{1}{|\Omega|}\int_\Omega u\,dx.
$$
The local attractor bifurcation and phase transition of the system have been extensively studied in Ma and Wang  \cite{MW2}. Other results relates to bifurcation of the problem can be found in \cite{BDW,Mis2}, etc.
 Here by applying the  theoretical results obtained above,  we try to provide some new results about  the dynamic bifurcation of the system and demonstrate   global features of the bifurcations.

\subsection{Mathematical setting of the system}

Denote by $(\.,\.)$ and $|\.|$ the usual inner product and norm of $L^2(\Omega)$, respectively.
For mathematical setting, we introduce the Hilbert space $H$  as follows:
$$
H=\{u\in L^2(\Omega):\,\,\,m(u)=0\}.$$ Let $A_0=-\De$ be the Laplacian in $H$  associated with the homogeneous boundary condition$$
\frac{\pa u}{\pa \nu}=0,\Hs x\in\pa\Omega.$$

Set $A=A_0^2$. Then $A$ is a positive-definite self-adjoint  operator in $H$ (and hence is a sectorial operator) with compact resolvent, and
$$\ba{ll}
D(A)=\left\{u\in H^4(\Omega)\cap H\mid\,\,\, \frac{\pa u}{\pa \nu}=\frac{\pa (\De u)}{\pa \nu}=0\mb{ on }\pa\Omega\right\}.\ea$$
The spectral $\sig(A_0)$ of  $A_0$ consists of a countably infinitely many  eigenvalues: $$0<\mu_1<\mu_2<\cdots< \mu_k\ra+\8.$$
Let $V:=D(A_0)=D(A^{1/2})$.
Denote  $||\.||$ the  norm in $V$. 

Define
$$
g_\lam(u)=\De(b_2 u^2+b_3 u^3), \Hs u\in V.$$
Then $g_\lam: V\ra H$ is locally Lipschitz, and the system (\ref{e:6.2}) can be reformulated in an abstract form:
\be\label{e:6.3}
u_t+ L_\lam u=g_\lam(u),\ee
 where $L_\lam=A_0^2-\lam A_0$.  We infer from Henry \cite{Henry}, Chap. 3 that for each $u_0\in V$, (\ref{e:6.3}) has a unique global strong solution $u(t)$ in $V$ with $u(0)=u_0$.

It is worth noticing that  the problem has a natural Lyapunov function $J(u)$,
$$
J(u)=\frac{1}{2}|\nab u|^2+\int_\Omega F_\lam(u)\,dx,
\hs \mb{
where }\,
F_\lam(s)=-\frac{\lam}{2}s^2+\frac{b_2}{3}s^3+\frac{b_3}{4}s^4.
$$

\subsection{Bifurcation from the trivial solution}

It is obvious that each eigenvector $w$ of $A_0$ corresponding to $\mu_k$ is also an eigenvector of $L_\lam$  corresponding to the  eigenvalue $$\b_k(\lam):=\mu_k^2-\lam\mu_k=\mu_k(\mu_k-\lam). $$ Because $H$ has a canonical basis consisting of eigenvectors of $A_0$, we deduce that
$
\b_k(\lam)$ $(k=1,2,\cdots$) are precisely all the eigenvalues  of $L_\lam$.

Let $\Phi_\lam$ be the semiflow generated by the system. We have

\bt\label{t:6.2} Assume  $b_2\ne 0$. Suppose  $A_0$ has an eigenvector $w$  corresponding to $\mu_j$ such that $\int_\Omega w^3\,dx\ne0$, and that
 $0$ is an isolated equilibrium of $\Phi_{\mu_j}$.

 Then
there exist a closed  neighborhood $U$ of $0$ in $V$ and a two-sided neighborhood $I_2$ of $\mu_j$  such that $\Phi_\lam$ has a nonempty maximal compact invariant set $K_\lam$ in $U\sm\{0\}$ for each $\lam\in I_2\sm\{\mu_j\}$.
Consequently for   $\lam\in I_2\sm\{\mu_j\}$, $\Phi_\lam$ has at least one nontrivial equilibrium.
\et
\noindent{\bf Proof.}  Since the system is a gradient-like one, by assumption   it is easy to check that $S_0=\{0\}$ is  an isolated invariant set of $\Phi_{\mu_j}$. In what follows we check  that $S_0$ is neither an attractor nor a repeller of the restriction $\Phi_{\mu_j}^c$ of $\Phi_{\mu_j}$ on $\cM^c$, and hence the conclusion of theorem immediately follows  from Theorem \ref{t:2.7}.

Denote $E_j$ the space spanned by the eigenvectors of $A_0$ corresponding to $\mu_j$. Then $H=E_j\bigoplus E_j^\bot$. Let $V_j^\bot=V\cap E_j^\bot$. Then  $V=E_j \bigoplus V_j^\bot$.
We infer from \cite{Ryba}, Chap. II, Theorem 2.1 that there is a small neighborhood $W$ of $0$ in $E_j$ and a $C^1$ mapping $\xi: W\ra V_j^\bot$ with
$$
\xi(v)=0(||v||^2)\hs(\mb{as }||v||\ra 0)
$$
such that $\cM^c=\{v+\xi(v):\,\,v\in W\}$ is a local center manifold of $\Phi_{\mu_j}$.

For $u=v+\xi(v)$, where $v\in W$, simple computations show that
$$
J(u)=\frac{1}{2}\int_\Omega |\xi'(v)\nab v|^2\,dx+\frac{b_2}{3}\int_\Omega v^3\,dx+o(||v||^3).
$$
 Here we have used the  facts that $\xi(v),\De \xi(v)\in E_j^\bot$. Setting $v=\tau w$, where $w$ is the eigenvector of $A_0$ given in the theorem, then
since $\xi'(v)=0(||v||)$, we have \be\label{e7.16}
 J(\tau w)=\tau^3\,\frac{b_2}{3}\int_\Omega w^3\,dx+o(|\tau|^3)\hs\mb{as }\tau\ra0.
\ee
As $\frac{b_2}{3} \int_\Omega w^3\,dx\ne 0$, by \eqref{e7.16} it is clear that $0$ is neither a local maximum  nor  minimum point of $J$, which completes the proof of what we desired.  \,$\Box$
\Vs
The following result demonstrate some global features of the dynamic bifurcation of the system.

\bt\label{t:6.3} Suppose $0$ is an isolated equilibrium of $\Phi_{\mu_j}$\,.   Let $\Gamma$ be the bifurcation branch in $V$ from the bifurcation point $(0,\mu_j)$. Set
$$
\Lam_0=\inf\{\lam:\,\,\Gamma_\lam\ne\emp\},\hs \Lam_1=\sup\{\lam:\,\,\Gamma_\lam\ne\emp\}\,,
$$
where $\Gamma_\lam=\{u:\,\,(u,\lam)\in\Gamma\}$ is the $\lam$-section of $\Gamma$.

Then $-\8<\Lam_0<\Lam_1\leq+\8$, and one of the  following assertions hold.

\begin{enumerate}
\item[(1)] $\Lam_1=+\8$.
\vs
\item[(2)] $ 0\in\ol{\Gamma_{\mu_j}\sm\{0\}}$.
\vs
\item[(3)]  There exists $\lam_1\ne \mu_j$ such that $(0,\lam_1)\in \Gamma$. Furthermore, either (i) there is a sequence $(u_k,\nu_k)\in\Gamma$ approaching $(0,\lam_1)$, where $u_k$ is a nontrivial equilibrium of $\Phi_{\nu_k}$ for each $k$;  or (ii) $\Gamma_{\lam_1}$ contains at least two distinct complete trajectories $\sig^\pm$ such that
$$J(\a(\sig^+))\equiv \mb{const.}>0, \hs \omega(\sig^+)=\{0\},$$
$$J(\omega(\sig^-))\equiv \mb{const.}<0, \hs \a(\sig^-)=\{0\}.$$
\end{enumerate}
\et

\br It is worth noticing that both $\a(\sig^+)$ and $\omega(\sig^-)$ in  (3) consist of  nontrivial equilibrium points. Therefore when (3) occurs,  $\Phi_{\lam_1}$ has at least two distinct nontrivial equilibria. When $\Gamma_{\lam_1}$ contains only a finite number of equilibria,  each of the two limit sets $\a(\sig^+)$ and $\omega(\sig^-)$ consists of exactly one equilibrium. Consequently $\sig^\pm$ become heteroclinic orbits.
\er
\noindent{\bf Proof of Theorem \ref{t:6.3}.} It can be easily shown that if $\lam<0$ is large enough, then the trivial solution $0$ is the global attractor of $\Phi_\lam$. Hence we necessarily have $\Lam_0>-\8$. The existence of local bifurcation branch also implies $\Lam_0<\Lam_1$.

Assume $\Lam_1<+\8$ (otherwise (1) holds true, and thus we are done).   Then $I=[\Lam_0,\Lam_1]$ is a compact interval.  Therefore we infer from the proof for the existence of a global attractor of the system  in Temam \cite{Tem} (see also \cite{LZ} etc.) that the system is dissipative uniformly with respect to $\lam\in I$. Specifically, there is a bounded set $B\subset V$ such that \be\label{e7.0}\cA_\lam\subset B,\Hs \A\,\lam\in I,\ee where $\cA_\lam$ is the global attractor of $\Phi_\lam$. Thus the bifurcation branch $\Gamma$ is bounded. Hence by Theorem \ref{gbt} we conclude that either (2) holds, or there is a $\lam_1\ne \mu_j$ such that $(0,\lam_1)\in \Gamma$.
To complete the proof of the theorem, there remains to check the alternatives in (3).

So we assume that  $(0,\lam_1)\in \Gamma$ for some $\lam_1\ne \mu_j$.
Suppose the first case (i) in (3) does not occur. Then  $0$ is  an isolated equilibrium of $\Phi_{\lam_1}$.
 Fix a $\de_1>0$ such that $\Phi_{\lam_1}$ has no equilibria other than the trivial one in the $\de_1$-neighborhood $\mB_{\de_1}$ of $0$ in $V$.
 By the definition of  bifurcation branch we deduce that  there exists a sequence $\nu_k\ra \lam_1$ such that for each  $k$, $\Phi_{\nu_k}$ has a nonempty compact invariant set $M_k\subset \Gamma_{\nu_k}$ with $0\not\in M_k$ such that
\be\label{e7.1}
\lim_{k\ra\8}d(0,M_k)=0.
\ee
For  convenience,  denote  $\sE(\Phi_\lam,M)$ the set of equilibria  of $\Phi_\lam$ in  $M\subset V$.
Let $$\sE_k:=\sE\(\Phi_{\nu_k},M_k\).$$ Then $\sE_k$ is a nonempty compact subset of $M_k$. As we have assumed that (i) does not occur, it can be easily seen that
 there exists $0<\de<\de_1$ such that
\be\label{e:6.10}
\liminf_{k\ra\8}d(0,\sE_k)\geq4\de>0.
\ee

  By \eqref{e7.1}, for  each $k$ we can  pick a $u_k\in M_k$ such that the sequence $u_k\ra0$ as $k\ra\8$. It can be assumed that \be\label{e7.2}||u_k||<\de\ee  for all $k$ (hence  $d(u_k,\sE_k)>3\de $).  Let $\gamma_k$ be a complete trajectory of $\Phi_{\nu_k}$ contained in $M_k$ with $\gamma_k(0)=u_k$. We have  \be\label{e:6.11}
\min_{t\leq0}J(\gamma_k(t))=J(\gam_k(0))=J(u_k)\ra0,\hs\mb{as }\,k\ra\8.
\ee
Set
$$
t_k=\min\{s<0:\,\,\max_{t\in[s,0]}||\gamma_k(t)-u_k||\leq 2\de\}.
$$
Noticing that $\a(\gamma_k)\subset \sE_k$, we deduce by (\ref{e:6.10}) and \eqref{e7.2}  that $t_k>-\8$,
and hence  $||\gamma_k(t_k)-u_k||=2\de.$ Thereby
\be\label{e7.13}\de\leq||\gamma_k(t_k)||\leq 3\de,\Hs k\geq 1.
\ee

Define a sequence of  complete trajectories $\sig_k$ as
\be\label{e:6.12}
\sig_k(t)=\gamma_k(t_k+t),\Hs t\in \R.
\ee
Since   all these trajectories are contained in the bounded set $B$ in \eqref{e7.0}, 
by  very standard argument  it can be shown that $\sig_k$ has a subsequence (still denoted by  $\sig_k$)  converging uniformly on any compact interval  to a complete trajectory $\sig^+$. It is trivial to check that  $\sig^+$ is contained in $\Gamma_{\lam_1}$. Observing that $\sig^+(0)=\lim_{k\ra\8}\gam_k(t_k)$, by \eqref{e7.13} we deduce that   \be\label{e7.15}\de\leq ||\sig^+(0)||\leq 3\de.\ee
Because
$$
J(\sig_k(0))\geq J(\sig_k(-t_k))=J(\gam_k(0))=J(u_k)\ra 0
$$
as $k\ra\8$, we also have  $J(\sig^+(0))\geq 0$.

On the other hand, as $\Phi_{\lam_1}$ has no equilibrium in $\mB_{\de_1}\sm \{0\}$,  by  \eqref{e7.15} we see  that $\sig^+(0)$ is not an equilibrium of $\Phi_{\lam_1}$.  Hence there is a small open  interval $I_\ve=(-\ve,\ve)$ such that 
 $J(\sig^+(t))$ is strictly decreasing in $t$ on $I_\ve$. Consequently $$J(\a(\sig^+))\equiv \mb{const.}>J(\sig^+(0))\geq 0.$$

In what follows we show that $\omega(\sig^+)=\{0\}$. If  $t_k$  has a bounded subsequence (still denoted by $t_k$) with $t_k\ra -\tau\leq 0$, then $$
 \sig^+(\tau)=\lim_{k\ra\8}\sig_k(-t_k)=\lim_{k\ra\8}\gam_k(0)=\lim_{k\ra\8}u_k=0.
$$
Hence $\sig^+(t)\equiv 0$ for $t\geq\tau$, which contradicts \eqref{e7.15}. Thus we know that  $t_k\ra-\8$. Since
$
||\gam_k(t)-u_k||\leq 2\de$ for $ t\in[t_k,0],
$
we have
$$||\gam_k(t)||\leq ||u_k||+ 2\de\leq 3\de,\Hs t\in[t_k,0].$$
Thereby
$$||\sig_k(t)||\leq 3\de,\Hs t\in[0,-t_k],$$
 from which it follows that $||\sig^+(t)||\leq 3\de$ for all $t\geq 0$. As $0$ is the unique equilibrium of $\Phi_{\lam_1}$ in $\mB_{\de_1}$ and $3\de<\de_1$, we immediately conclude that
$\omega(\sig^+)=\{0\}$.

\vs Likewise, we can prove that there is a complete trajectory  $\sig^-$ in $\Gamma_{\lam_1}$ such that
$$J(\omega(\sig^-))\equiv \mb{const.}<0, \hs \a(\sig^-)=\{0\}.$$
The proof of the theorem is finished. $\Box$

\br We have assumed in Theorems \ref{t:6.2} and  \ref{t:6.3} that the trivial  solution $0$ of the system  is an isolated equilibrium of $\Phi_{\mu_j}$. In general it seems to be difficult to verify this condition due to the degeneracy. However, in some particular but important cases one can really do so. For instance, if $b_2=0$ then it can be shown that the equilibrium $0$  is  isolated with respect to $\Phi_{\mu_j}$ (see the proof of Theorem 9.4 in Ma and Wang \cite{MW1b}).
\er

\vskip 0.2in

\noindent{\bf Acknowledgements:} D. Li is supported by NSFC-10771159,
11071185 and Z.-Q. Wang is partially supported by NSFC-11271201.
The authors would like to express their gratitude to the referees for their valuable comments and suggestions which helped  them greatly improve the quality  of the paper.

\end{document}